\documentclass[11pt]{article}
\usepackage{amsmath, graphicx, amsfonts,amssymb, calrsfs}
\usepackage{amsfonts,mathrsfs, color, amsthm}
\usepackage{enumitem}
\usepackage{hyperref}
\hypersetup{
  colorlinks   = true, 
  urlcolor     = blue, 
  linkcolor    = blue, 
  citecolor   = red 
  }

\usepackage{adjustbox}

\newcommand{\vleq}{\mathrel{\rotatebox{90}{$\geq$}}}

\newtheorem{theorem}{Theorem}[section]
\newtheorem{lemma}{Lemma}[section]
\newtheorem{remark}{Remark}[section]
\newtheorem{proposition}{Proposition}[section]
\newtheorem{corollary}{Corollary}[section]

\newtheorem{definition}{Definition}[section]

\def\bt{\begin{theorem}}
\def\et{\end{theorem}}
\def\bl{\begin{lemma}}
\def\el{\end{lemma}}
\def\br{\begin{remark}}
\def\er{\end{remark}}
\def\bc{\begin{corollary}}
\def\ec{\end{corollary}}
\def\bd{\begin{definition}}
\def\ed{\end{definition}}
\def\bp{\begin{proposition}}
\def\ep{\end{proposition}}

\def\bpf{\begin{proof}}
\def\epf{\end{proof}}
\def\be{\begin{equation}}
\def\ee{\end{equation}}
\def\bea{\begin{eqnarray}}
\def\eea{\end{eqnarray}}

\def\cK{\mathcal{K}}
\def\sphere{S^{n-1}}

\def\Rn{{\mathbb R^n}}
 \def\R{\mathbb{R}}
\def\cH{\mathcal{H}}

\def\affp{C_{p, \tau}}

\begin{document}
\title{Sharp geometric inequalities for the general $p$-affine capacity \footnote{Keywords: Asymmetric $L_p$ affine Sobolev inequality, general $L_p$ affine isoperimetric inequality,  isocapacitary inequality, $L_p$ affine isoperimetric inequality, $L_p$ affine Sobolev inequality,  $L_p$ projection body,  $p$-affine capacity,  $p$-integral affine surface area, $p$-variational capacity.}}
\author{Han Hong and Deping Ye}
\date{}
\maketitle

\begin{abstract} 
In this article, we propose the notion of the general $p$-affine capacity and prove some basic properties for the general $p$-affine capacity, such as affine invariance and monotonicity. The newly proposed general $p$-affine capacity is compared with several classical geometric quantities, e.g., the volume, the $p$-variational capacity and the $p$-integral affine surface area. Consequently, several sharp geometric inequalities for the general $p$-affine capacity are obtained. These inequalities extend and strengthen many well-known  (affine) isoperimetric and  (affine) isocapacitary inequalities.  

\vskip 2mm 2010 Mathematics Subject Classification: 46E30, 46E35, 52A38, 53A15.
\end{abstract}

\section{Introduction}

Many objects of interest and fundamental results in convex geometry are related to the $L_p$ projection bodies \cite{LXZ, LYZ1, LYZ-duke, Petty1971}.  For $p\geq 1$, the $L_p$ projection body of a convex body (i.e., a compact convex subset with nonempty interior) $K\subset \Rn$ containing the origin in its interior is determined by its support function $h_{\Pi_{p}(K)}: \sphere\rightarrow \R$,  whose definition is formulated as follows (up to a multiplicative constant):  for any $\theta \in \sphere$,  \begin{eqnarray*} h_{\Pi_{p}(K)}(\theta)= \bigg( \int_{\partial K}\Big(\frac{\theta \cdot \nu_K(x)}{2}\Big)^p\cdot |x\cdot \nu_K(x)|^{1-p}\, d\cH^{n-1}(x)\bigg)^{\frac{1}{p}}\end{eqnarray*}  with $\nu_K$ the unit outer normal vector of $K$ at $x\in \partial K$  and $\cH^{n-1}$ denotes the $(n-1)$-dimensional Hausdorff measure of $\partial K$, the boundary of $K$  (see Section \ref{2-background-1} for details on the notations).  Define $\Phi_p(K)$, the $p$-integral affine surface area of $K$,  by $$\Phi_p(K)= \left(\int_{S^{n-1}} \big[ h_{\Pi_{p}(K)}(u)\big]^{-n} \, du\right)^{-\frac{p}{n}}$$ where $\,du$ is the normalized spherical measure on the unit sphere $\sphere$. Let $B_n$ be the unit Euclidean ball in $\Rn$ and $V(K)$ denote the volume of $K$. The following $L_p$ affine isoperimetric inequality for the $p$-integral affine surface area holds \cite{LXZ, LYZ1, LYZ-duke, Petty1971,  Zhang1}: for $p\geq 1$ and for $K$ a convex body with the origin in its interior,  \begin{equation}\label{classical-Lp-affine-isoperimetric}  \bigg(\frac{\Phi_{p}(K)}{\Phi_{p}(B_n)}\bigg)^{\frac{1}{n-p}}\geq   \bigg(\frac{V(K)}{V(B_n)}\bigg)^{\frac{1}{n}}\end{equation}   with equality if and only if $K=TB_n$ if $p>1$  and $K=TB_n+x_0$ if $p=1$ for some invertible linear transform $T$ on $\Rn$  and some $x_0\in \Rn$. Note that inequality (\ref{classical-Lp-affine-isoperimetric}) is invariant under the volume preserving linear transforms and hence is stronger than the well-known $L_p$ isoperimetric inequality \cite{Gardner-1, Lutwak, Schneider}:  \begin{equation}\label{classical-Lp-isoperimetric}  \bigg(\frac{S_{p}(K)}{S_{p}(B_n)}\bigg)^{\frac{1}{n-p}}\geq   \bigg(\frac{V(K)}{V(B_n)}\bigg)^{\frac{1}{n}}\end{equation}  with equality if and only if $K$ is an Euclidean ball in $\Rn$ (if $p>1$, the center needs to be at the origin). Here $S_p(K)$ is the $p$-surface area of $K$ and can be formulated by  \begin{equation} \label{definition-p-surface area} S_p(K)=\int_{\partial K} |x \cdot \nu_K(x)|^{1-p}\,d\cH^{n-1}(x).\end{equation} 

It is well known that inequality (\ref{classical-Lp-isoperimetric}) can be strengthened by the isocapacitary inequality related to the $p$-variational capacity. For a compact set $K\subset \Rn$, its  $p$-variational capacity, denoted by  $C_p(K)$, can be formulated by   (see e.g. \cite{Evan, Mazya-85, Mazya})  \begin{eqnarray*} C_p(K) =\inf \Big\{ \int_{\R^n} |\nabla f|^p\,dx: \  \ f\in C_c^{\infty} \ \text{and}\ f\geq  1  \ \text{on} \ K \Big\},\end{eqnarray*}   where $\nabla f$ denotes the gradient of $f$ and $C^{\infty}_{c}$ is the set of smooth functions with compact supports in $\R^n$. The $p$-variational capacity is an important geometric invariant which has close connection with the $p$-Laplacian partial differential equation and has important applications in many areas, e.g., analysis, mathematical physics and partial differential equations (see e.g., \cite{Evan, Mazya-85, Mazya} and references therein).  In particular, the Brunn-Minkowski type inequalities and the Hadamard variational formulas for the $p$-variational capacity have been established in, e.g., \cite{Borell, Caffarelli1996, CNSXYZ, CS2003, HYZ-17, Jerison, Jerison-1996, xiongzou}.  The following inequality for the $p$-variational capacity holds \cite{LXZ, Mazya-85}:  for $p\in [1, n)$ and for $K$ being a Lipschitz star body with the origin in its interior, \begin{equation}\label{isocapacity-11-introduction}   \bigg(\frac{S_{p}(K)}{S_{p}(B_n)}\bigg)^{\frac{1}{n-p}}\geq  \bigg(\frac{C_p(K)}{C_p(B_n)}\bigg)^{\frac{1}{n-p}}\geq  \bigg(\frac{V(K)}{V(B_n)}\bigg)^{\frac{1}{n}}. \end{equation}  The $p$-variational capacity behaves rather similar to the $p$-surface area and is lack of the affine invariance. Very recently, Xiao \cite{Xiao2015, Xiao} introduced  an affine relative of the $p$-variational capacity and named it as the $p$-affine capacity. This new notion is denoted by $C_{p, 0}(K)$ in this article and its definition is equivalent to, as proved in Section \ref{section:3-definition}, the following: for $p\in [1, n)$ and for $K$ a compact set in $\Rn$,   $$
C_{p, 0}(K)=\inf \Big \{\mathcal{H}_p(f): \ f\in C_c^{\infty} \ \text{and}\ f\geq  1 \ \ \text{on} \ \ K  \Big\}, 
$$ where $\mathcal{H}_p(f)$ is the $p$-affine energy of $f$:  $$\mathcal{H}_p(f)= \bigg(\int_{S^{n-1}}\bigg(\int_{\R^n} \frac{|\langle u\cdot \nabla f\rangle|^p}{2} dx\bigg)^{-\frac{n}{p}}du\bigg)^{-\frac{p}{n}}.$$ The following affine isocapacity inequality was also established in \cite[Theorems 3.2 and  3.5]{Xiao} and \cite[Theorems 1.3' and 1.4']{Xiao2015}:   for $p\in [1, n)$ and for $K$ an origin-symmetric convex body, one has \begin{equation}\label{affine-isocapacity-11-introduction} \bigg(\frac{\Phi_{p}(K)}{\Phi_{p}(B_n)}\bigg)^{\frac{1}{n-p}}\geq \bigg(\frac{C_{p, 0}(K)}{C_{p, 0}(B_n)}\bigg)^{\frac{1}{n-p}}\geq \bigg(\frac{V(K)}{V(B_n)}\bigg)^{\frac{1}{n}}. \end{equation}  The second inequality of (\ref{affine-isocapacity-11-introduction}) indeed also holds for any compact set $K\subset \Rn$.   Again inequality (\ref{affine-isocapacity-11-introduction}) is invariant under the volume preserving linear transforms and hence is stronger than inequality (\ref{isocapacity-11-introduction}).  Moreover, inequality (\ref{affine-isocapacity-11-introduction}) can be viewed as the affine relative of  inequality (\ref{isocapacity-11-introduction}).  See e.g., \cite{WangTuo2017, Xiao-1, Xiaozhang} for more works related to affine capacities. We would like to mention that the $p$-affine energy is the key ingredient in many fundamental analytical inequalities, see e.g., \cite{CLYZ, HJM2016, LYZ, VHN, WangTuo2012, WangTuo2013, Xiao2007,Zhang}. 
  
  It is our goal in this article to study a concept more general than the $p$-affine capacity and to establish stronger sharp geometric inequalities. The motivation is a result from recent studies, such as,  the general $L_p$ affine isoperimetric inequalities and asymmetric affine $L_p$ Sobolev inequalities by Haberl and Schuster \cite{HSgeneral, HSasymmetric},   asymmetric affine P\'{o}lya-Szeg\"{o} principle by Haberl, Schuster and Xiao \cite{HSXasymmetric}  and Minkowski valuations by Ludwig \cite{Ldvaluation}.  The key in \cite{HSgeneral} is to replace  $h_{\Pi_{p}(K)}$ by its asymmetric counterpart $h_{\Pi_{p, \tau}(K)}: \sphere \rightarrow \R$: for any $p\geq 1$, for any $\tau\in [-1, 1]$ and for $K$ a convex body with the origin in its interior,  $$\big[h_{\Pi_{p, \tau}(K)}(\theta)\big]^p= \int_{\partial K} \big[\varphi_\tau(\theta\cdot \nu_K(x))\big]^p \cdot |x\cdot \nu_K(x)|^{1-p}\, d\cH^{n-1}(x)$$ for $\theta \in \sphere$,  where  \begin{equation}\label{function-tau-p}
\big[\varphi_\tau(t)\big]^p=\Big(\frac{1+\tau}{2}\Big) t_+^p +\Big(\frac{1-\tau}{2}\Big) t_-^p
\end{equation} with  $t_+=\max\{0,  t\}$ and  $t_-=\max\{0,  -t\}$ for any $t\in \R$. We point out that this extension is a key step from the $L_p$ Brunn-Minkowski theory of convex bodies to the Orlicz theory and its dual (see e.g., \cite{GHW2014, ghwy15,  Ludwig2010, LYZ2010a, LYZ2010b, XJL, Zhub2014}).  Similarly, the key in \cite{HSasymmetric, HSXasymmetric}  is to replace the $p$-affine energy function $\mathcal{H}_p(f)$ by its asymmetric counterpart: for any $p\in [1, n)$, for any $\tau\in [-1, 1]$ and for any $f\in C_c^{\infty}$, $$ \mathcal{H}_{p,\tau}(f) =  \left(\int_{S^{n-1}} \left(\int_{\R^n} \big[\varphi_\tau(\nabla_uf)
\big]^p \, dx\right)^{-\frac{n}{p}}\,du\right)^{-\frac{p}{n}}. 
$$   When $\tau=0$,  $\mathcal{H}_{p,\tau}(f)$ goes back to the $p$-affine energy $\mathcal{H}_p(f)$. It is worth to mention that to deal with $\mathcal{H}_{p,\tau}(f)$ is much more challenging than $\mathcal{H}_p(f)$, mainly because the $L_p$ convexifications of level sets of a smooth function $f$ in the latter case always contain the origin in their interiors but in the former may not contain the origin in their interiors.  These asymmetric extensions have also been widely used  to study affine Sobolev type inequalities, the affine P\'{o}lya-Szeg\"{o} principle  as well as many other affine isoperimetric inequalities, see  e.g., \cite{VHN,Ober, WangTuo2015, Weberndorfer}.

In Section \ref{section:3-definition}, we provide several equivalent definitions for  the general $p$-affine capacity, which will be denoted by $\affp(\cdot)$.  One of them reads:  for any $p\in [1, n)$, for any $\tau \in [-1, 1]$ and for any compact set $K\subset \Rn$, $$
\affp(K)=\inf \Big \{\mathcal{H}_{p,\tau}(f): \ f\in C^{\infty}_c \ \text{and}\ f\geq 1 \ \text{on} \ K  \Big\}.$$ Basic properties for the general $p$-affine capacity, such as, monotonicity,  affine invariance, translation invariance, homogeneity and the continuity from above, are established in Section \ref{section:4}. Similarly, the general $p$-integral affine surface area of a Lipschitz star body $K$ is defined in Subsection  \ref{section:4-1} by: for any $p\in [1, n)$ and for any $\tau \in [-1, 1]$,  $$\Phi_{p, \tau} (K)= \left(\int_{S^{n-1}} \big[h_{\Pi_{p, \tau}(K)}(u)\big]^{-n} \, du\right)^{-\frac{p}{n}}.$$ Note that when $\tau=0$, then $\Phi_{p, 0}(K)=\Phi_p(K).$  The sharp geometric inequalities for the general $p$-affine capacity are established in Section \ref{section:5}. Roughly speaking,  for $K$ a convex body containing the origin in its interior, these sharp geometric inequalities can be summarized as follows: for all $p\in [1, n)$ and for all $0\leq \tau \leq \eta \leq 1$, then  
  \begin{equation}\label{chain-inequalities}  \begin{array} {cccc}  
 \Big(\frac{V(K)}{V(B_n)}\Big)^{\frac{1}{n}} &\leq \Big(\frac{C_{p, \eta}(K)}{C_{p, \eta}(B_n)}\Big)^{\frac{1}{n-p}} &\leq \Big(\frac{\Phi_{p, \eta}(K)}{\Phi_{p, \eta}(B_n)}\Big)^{\frac{1}{n-p}} & \leq \Big(\frac{S_{p}(K)}{S_{p}(B_n)}\Big)^{\frac{1}{n-p}}  
 \\  [10pt]  &\vleq &\vleq &  
 \\ \Big(\frac{V(K)}{V(B_n)}\Big)^{\frac{1}{n}} &\leq \Big(\frac{C_{p, \tau}(K)}{C_{p, \tau}(B_n)}\Big)^{\frac{1}{n-p}} &\leq \Big(\frac{\Phi_{p, \tau}(K)}{\Phi_{p, \tau}(B_n)}\Big)^{\frac{1}{n-p}} &\leq \Big(\frac{S_{p}(K)}{S_{p}(B_n)}\Big)^{\frac{1}{n-p}}.
 \end{array} \end{equation} Inequality (\ref{affine-isocapacity-11-introduction}) turns out to be a special (and indeed the maximal) case of the above chain of inequalities. Hence,  (\ref{chain-inequalities}) extends and strengthens many well-known  (affine) isoperimetric and (affine) isocapacitary inequalities, such as,    \cite[Theorem 1]{HSgeneral} by Haberl and Schuster,   \cite[inequality (13)]{LXZ} by Ludwig, Xiao and Zhang, and \cite[Theorems 3.2 and  3.5]{Xiao} by Xiao. Moreover, we also prove that, for any $p\in [1, n)$ and for any $\tau\in [-1, 1]$,  \begin{equation} \label{affine-capacity vs capacity} \bigg(\frac{S_p(K)}{S_p(B_n)}\bigg)^{\frac{1}{n-p}}\geq \bigg(\frac{C_p(K)}{C_p(B_n)}\bigg)^{\frac{1}{n-p}}\geq \bigg(\frac{\affp(K)}{\affp(B_n)}\bigg)^{\frac{1}{n-p}}\geq \bigg(\frac{V(K)}{V(B_n)}\bigg)^{\frac{1}{n}}, \end{equation} which extends and strengthens, e.g., inequality (\ref{isocapacity-11-introduction}),  \cite[(12)]{LXZ} by Ludwig, Xiao and Zhang, and \cite[Remark 2.7]{Xiao} by Xiao.  Note that inequalities (\ref{chain-inequalities}) and  (\ref{affine-capacity vs capacity}) work for more general compact sets than convex bodies, and we will explain the details in Section  \ref{section:5}.

\section{Background and Notations}\label{2-background-1}

A compact set $M\subset \Rn$ is said to be a star body (with respect to the origin $o$) if the line segment jointing $o$ and $x$, for all $x\in M$, is contained in $M$. For each star body $M$, one can define the radial function $\rho_M$ of $M$ as follows: for all $x\in \Rn\setminus \{o\}$, $$\rho_M(x)=\mathrm{max}\{\lambda\geq 0:\ \lambda x\in M \}. $$ The star body $M$ is said to be a Lipschitz star body if the boundary of $M$ is Lipschitz.

A compact convex subset in $\R^n$ with nonempty interior is called a convex body. By $\cK_0$,  we mean the set of all convex bodies with the origin $o$ in their interiors. Each $K\subset \cK_0$ is (uniquely) associated with two 
continuous functions defined on the unit sphere $\sphere$: the radial function $\rho_K$ and  the support function $h_K$. Hereafter,  for $u\in \sphere$,  \begin{eqnarray}
	h_K(u)&=&\mathrm{max}\{y\cdot u: \ y\in K\}, \nonumber
\end{eqnarray} where $x\cdot y$ is  the standard inner product of $x$ and $y$ in $\Rn$.  The support function  $h_K: \sphere\rightarrow (0, \infty)$ of a convex body $K\in \cK_0$ can be extended to $\Rn\setminus\{o\}$ as follows: $h_K(x)=rh_K(u)$ for any $x\in \Rn \setminus\{o\}$ with $x=ru$. It can be easily checked that the extended function $h_K: \Rn\setminus\{o\}\rightarrow (0, \infty)$  is sublinear, i.e., $h_K$ has the positive homogeneity of degree 1 and satisfies $$h_K(x+y)\leq h_K(x)+h_K(y)$$ for all $x, y\in \Rn\setminus\{o\}$. On the other hand, if a function $h: \Rn\setminus \{o\} \rightarrow (0,\infty)$ is sublinear, then $h$ is the support function of a convex body $K\in \cK_0$ \cite{Schneider}.  For each $K\in \cK_0$, its polar body $K^\circ$ is $$
K^\circ=\{x\in \R^n: \ x\cdot y \leq 1  \ \ \mathrm{for\ \ all} \  \  y\in K \}.
$$ It is easily checked that  \begin{equation}\label{supportradial}
\rho_{K^\circ}=\frac{1}{h_K} \ \ \ \ \mathrm{and} \ \ h_{K^\circ}=\frac{1}{\rho_K} .
\end{equation}

The standard notation $\cH^{k}$ is for the $k$-dimensional Hausdorff measure. In the case of $k=n$, we use $V(\cdot)$ to denote the volume instead of $\cH^n$. In particular, the volume of  the unit Euclidean ball $B_n$, denoted by $\omega_n$ for simplicity, has the following expression: $$\omega_n=\frac{\pi^{n/2}}{\Gamma(1+n/2)},$$ where $\Gamma(\cdot)$ is the Gamma function $$
\Gamma(x)=\int_{0}^{\infty}t^{x-1}e^{-t}\ dt.$$  The Beta function $B(\cdot, \cdot)$ is closely related to the Gamma function, and it has the form
$$  B(x,y)=\int_{0}^{1}t^{x-1}(1-t)^{y-1}\ dt.
$$ It is easily checked that 
\begin{equation*} 
B(x,y)=\frac{\Gamma(x)\Gamma(y)}{\Gamma(x+y)}.
\end{equation*}

It is convention to use $\,d\sigma$ for the spherical measure of $\sphere$. In later context, the normalized spherical measure $\,du$ is often used, i.e.,  $$ \,du=\frac{\,d\sigma}{n\omega_n} \ \ \ \text{and} \ \ \ \int_{S^{n-1}}du=1.$$  The volume of each $K\in \cK_0$ can be calculated by \begin{equation} \label{volumeformula} V(K)=\frac{1}{n}\int_{\sphere} \rho_K^n(u)\,d\sigma(u)  \ \ \ \mathrm{or}\ \ \ V(K)=\frac{1}{n}\int_{\sphere} h_K(u)\,dS(K, u), \end{equation} where $S(K,\cdot)$ is the classical surface area measure of $K\in \cK_0$ defined on $\sphere$. Denote by  $C(S^{n-1})$ the set of continuous functions on $\sphere$. The classical surface area measure    $S(K,\cdot)$ has the following analytic interpretation: for all $f\in C({S^{n-1}})$,   \begin{equation}\label{surfaceareameasure}
\int_{S^{n-1}} f(u) \ dS(K,u)=\int_{\partial K}f(\nu_K(x)) \ d\mathcal{H}^{n-1}(x),
\end{equation} where $\nu_K(x)$ is an outer unit normal vector  at $x\in \partial K$, the boundary of $K$.  For each $K\in \cK_0$,  $\nu_K(x)$ exists almost everywhere on $\partial K$ with respect to  $\mathcal{H}^{n-1}$ \cite{Schneider}.

A smooth function is a real valued function $f:\Rn\rightarrow \R$ which is infinitely continuously differentiable. Denote by $C^{\infty}$ the set of smooth functions with continuous derivatives of all orders, and by $C^{\infty}_{c}$ (or $C^{\infty}_c(\R^n)$)  the set of functions in $C^{\infty}$ with compact support in $\R^n$. The gradient of $f\in C^{\infty}_{c}$ is denoted by $\nabla f$.  For $1\leq p<\infty$ and   $f\in C^{\infty}_{c}$, consider the norm $$\|f\|_{1, p}=\|f\|_p+\|\nabla f\|_p=\bigg(\int_{\Rn} |f|^p\,dx\bigg)^{1/p}+\bigg(\int_{\Rn} |\nabla f |^p\,dx\bigg)^{1/p}.$$  We also use $\lVert f\lVert_\infty$  to denote the maximal value  (or supremum) of $|f|$. The closure of $C^{\infty}_{c}$ under the norm $\|\cdot\|_{1, p}$ is denoted by  $W_0^{1,p}$. Note that the Sobolev space $W_0^{1,p}$ is a Banach space and each $f\in W_0^{1,p}$ is a real valued $L_p$ function on $\R^n$ with weak $L_p$ partial derivative (see e.g. \cite{Evan} for more details about the Sobolev space). Hereafter, when $f\in W_0^{1,p}$ is not smooth enough, $\nabla f$ means the weak partial gradient.  By $\nabla_zf$ we mean the inner product of $z$ and $\nabla f$, namely  
$\nabla_zf=z\cdot \nabla f.$ When $u\in \sphere$, $\nabla_uf$ is just the directional derivative of $f$ along the direction $u$. Clearly $\nabla_zf$ is linear about $z\in \Rn$.

For a subset $E\subset \Rn$, $\mathbf{1}_E$ denotes the indicator function of $E$, that is, $\mathbf{1}_E(x)=1$ if $x\in E$ and otherwise $0$. Let $|x|=\sqrt{x\cdot x}$ be the Euclidean norm of $x\in \Rn$. The distance from a point $x\in \R^n$ to a subset $E\subset \R^n$, denoted by $\mathrm{dist}(x, E)$, is defined by
$$
\mathrm{dist}(x, E)=\mathrm{inf}\{|x-y|:\ y\in E \}.
$$
Note that if $x\in \bar{E}$, the closure of $E$,  then $\mathrm{dist}(x, E)=0$.

For any real number $t>0$, define the level set $[f]_t$ of $f\in C_{c}^{\infty}$ by
\begin{equation}\label{levelset}
[f]_t=\{x\in \R^n: |f(x)|\geq t\}.
\end{equation}  For all $t\in (0, \|f\|_\infty)$, $[f]_t$ is a compact set. The Sard's theorem implies that, for almost every  $t\in (0, \|f\|_\infty)$,  the smooth $(n-1)$ submanifold  $$\partial[f]_t=\{x\in \R^n: |f(x)|=t\}$$ has nonzero normal vector $\nabla f(x)$ for all $x\in \partial[f]_t$. Denoted by $\nu(x)=-\nabla f(x)/|\nabla f(x)|$ and \begin{equation*}
\{\nu(x): x\in \partial[f]_t\}=S^{n-1}.
\end{equation*}  An often used formula in our proofs is the well-known Federer's coarea formula (see \cite{Federer}, p.289): suppose that $\Omega$ is an open set in $\R$ and $f: \Rn\rightarrow \R$ is a Lipschitz function, then \begin{equation}\label{Federer's coarea formula}
\int_{f^{-1}(\Omega)\bigcap\{|\nabla f|>0\}}g(x)\ dx=\int_{\Omega}\int_{f^{-1}(t)}\frac{g(x)}{|\nabla f(x)|}\ d\mathcal{H}^{n-1}(x)\ dt,
\end{equation}
for any measurable function $g: \R^n\rightarrow[0,\infty)$.

Denote by $\R^*$ the subset of $\R$ that contains nonnegative real numbers.  Let
$\varphi_\tau: \R\rightarrow \R^*$ be  the function given by formula (\ref{function-tau-p}), that is, for $\tau\in[-1,1]$ and $t\in \R$,
 \begin{equation}\label{function-tau}
\varphi_\tau(t) =\Big(\frac{1+\tau}{2}\Big)^{1/p}  t_+  +\Big(\frac{1-\tau}{2}\Big)^{1/p}  t_-.
\end{equation}  It is easily checked that  $\varphi_\tau$ has positive homogeneous of degree 1 and subadditive, i.e.
\begin{equation}\label{taufunction}
\varphi_\tau(\lambda t)=\lambda\varphi_\tau(t) \ \ \ \mathrm{for} \ \ \ \lambda\geq 0 \ \ \mathrm{and} \ \ \ \varphi_\tau(t_1+t_2)\leq\varphi_\tau(t_1)+\varphi_\tau(t_1).
\end{equation}  Special cases, which are commonly used,  are $\varphi_0(t)=2^{-1/p}|t|$,  $\varphi_1(t)=t_+$ and  $\varphi_{-1}(t)=t_-$. We would like to mention that the function $\psi_{\eta}: \R\rightarrow \R^*$ for each $\eta\in [-1, 1]$ given by $$\psi_{\eta}(t)=|t|+\eta t$$ is also commonly used  in convex geometry (see e.g. \cite{HSgeneral, Ldvaluation}). However,  if we let $$\tau=\frac{(1+\eta)^p-(1-\eta)^p}{(1+\eta)^p+(1-\eta)^p},$$  then  $\psi_{\eta}^p=\big((1+\eta)^p+(1-\eta)^p\big)\cdot   \varphi_{\tau} ^p.$ In later context,  the theory for the general $p$-affine capacity will be developed only based on $\varphi_{\tau}$ because it is more convenient to prove the convexity or concavity of the general $p$-affine capacity with $\varphi_{\tau}$.  

  We shall need the following result  (see, e.g., \cite[Lemma 1.3.1 (ii)]{Groemer}), which is crucial in the computation of involved integral on $\sphere$.

\bl \label{calculation}
If  $v\in S^{n-1}$ and $\Phi$ is a bounded Lebesgue integrable function on $[-1,1]$, then $\Phi(u\cdot v)$, considered as a function of $u\in S^{n-1}$, is integrable with respect to the normalized spherical measure $\,du$. Moreover, 
$$
\int_{S^{n-1}} \Phi(u\cdot v) \, du=\frac{(n-1)\omega_{n-1}}{n\omega_n}\int_{-1}^{1}\Phi(t)(1-t^2)^{\frac{n-3}{2}}\, dt.
$$
\el 
 
It can be easily checked that for $p>0$ \begin{eqnarray*} 
	\int_{-1}^{1}  t_+^p  (1-t^2)^{\frac{n-3}{2}}\, dt  &=& \int_{-1}^{1}  t_-^p  (1-t^2)^{\frac{n-3}{2}}\, dt\\ &=&\int_{0}^{1} t^p  (1-t^2)^{\frac{n-3}{2}}\, dt\\ &=&\frac{1}{2}\cdot  \int_{0}^{1}t^{\frac{p+1}{2}-1} (1-t)^{\frac{n-1}{2}-1}\, dt\\ &=&\frac{1}{2}\cdot  B\Big(\frac{p+1}{2},  \frac{n-1}{2}\Big). 
\end{eqnarray*} In particular, if $\Phi=\varphi_\tau^p$,  it follows from (\ref{function-tau-p}) and Lemma \ref{calculation}  that, for $p> 0$ and for any $u\in \sphere$,   \begin{eqnarray}
\int_{S^{n-1}}[ \varphi_\tau(u\cdot v)]^p\, du    &=&\frac{(n-1)\omega_{n-1}}{n\omega_n}\int_{-1}^{1}\Big[\Big(\frac{1+\tau}{2}\Big) t_+^p +\Big(\frac{1-\tau}{2}\Big) t_-^p\Big] (1-t^2)^{\frac{n-3}{2}}\, dt \nonumber  \\ &=& \frac{(n-1)\omega_{n-1}}{2 n\omega_n}\cdot B\Big(\frac{p+1}{2},  \frac{n-1}{2}\Big)   \ \  \ \ \ (:= A(n, p)). 
\label{computation}
\end{eqnarray}

\section{The general $p$-affine capacity}\label{section:3-definition}

In this section, the general $p$-affine capacity is proposed and several equivalent formulas for  the general $p$-affine capacity are provided.  Throughout, the general $p$-affine capacity of a compact set $K\subset \Rn$ will be denoted by $\affp(K).$ For convenience, let $$\mathcal{E}(K)=\{f: f\in W_0^{1, p},\  f\geq \mathbf{1}_K \}.$$
For each $f\in W_0^{1, p}$, let  $\nabla^+_uf(x)=\max\{\nabla _u f(x), 0\}$, $\nabla^-_uf(x)=\max\{-\nabla _u f(x), 0\},$ and \begin{eqnarray} \mathcal{H}_{p,\tau}(f) &=& \left(\int_{S^{n-1}} \|\varphi_\tau(\nabla_uf)\|_p^{-n}\,du\right)^{-\frac{p}{n}} \nonumber  \\  &=&  \left(\int_{S^{n-1}} \left(\int_{\R^n} \Big[\Big(\frac{1+\tau}{2}\Big) (\nabla^+_uf(x))^p +\Big(\frac{1-\tau}{2}\Big) (\nabla^-_uf(x))^p\Big]\ dx\right)^{-\frac{n}{p}}\,du\right)^{-\frac{p}{n}}.  \ \ \ \ \ \ \ \ \label{definition-h-f}  
\end{eqnarray} 

\bd\label{definition1}
Let $K$ be a compact subset in $\R^n$ and the function $\varphi_{\tau}$ be as in (\ref{function-tau}). For $1\leq p<n$, define the general $p$-affine capacity of $K$ by
\begin{eqnarray*} 
	\affp(K)&=&\inf_{f\in \mathcal{E}(K)} \mathcal{H}_{p,\tau}(f).
\end{eqnarray*} \ed 

  \noindent {\bf Remark.}  
For any compact set $K\subset \Rn$ and for any $\tau\in [-1, 1]$, $\affp(K)<\infty$ if $p\in [1, n)$.  According to  the proofs of (\ref{monotone-p-affine-capacity}) and Theorem \ref{property}, the desired boundedness argument follows if  $\affp(B_n)<\infty$ is verified. To  this end, let  $K=B_n$  and  $\varepsilon>0$. Consider   
$$
f_\varepsilon(x)=\left\{\begin{array}{ll} 
0, &\ \mathrm{if}\  |x|\geq 1+\varepsilon, \\
1-\frac{|x|-1}{\varepsilon}, &\ \mathrm{if}\ 1<|x|<1+\varepsilon,\\ 1, &\  \mathrm{if}\ |x|\leq 1.  \end{array}\right. $$ It can be checked that  $f_{\varepsilon}\in W_0^{1, p}$ and $f_\varepsilon$ has its weak derivative to be 
$$
\nabla f_\varepsilon(x)= \left\{\begin{array}{ll} 
0, & \ \mathrm{if}\  |x|\notin (1,  1+\varepsilon), \\
-\frac{x}{\varepsilon |x|}, & \ \mathrm{if}\ |x|\in (1, 1+\varepsilon).  \end{array}\right.
$$  
This further implies that, together with Fubini's theorem, (\ref{taufunction})  and (\ref{computation}),  \begin{eqnarray*} \|\varphi_\tau(\nabla_uf_\varepsilon)\|_p^p &=&\int_{\R^n} \big[ \varphi_\tau(\nabla_uf_\varepsilon(x))\big]^p\, dx  \\ &=&\int_{\{x\in \Rn: 1< |x|<1+\varepsilon\}} \bigg[ \varphi_\tau\bigg(\! -\frac{u\cdot x}{\varepsilon |x|} \bigg)\bigg]^p\, dx  \\ &=&  \varepsilon^{-p} \int_1^{1+\varepsilon} r^{n-1}\,dr \cdot \int _{\sphere}  \big[\varphi_{\tau}( -u\cdot v) \big]^p\, d\sigma(v)\\ &=&  \frac{(1+\varepsilon)^n-1}{\varepsilon^{p}} \cdot  \omega_n\cdot  A(n, p).  \end{eqnarray*} 
It follows from (\ref{definition-h-f}) that  $$\mathcal{H}_{p,\tau}(f_\varepsilon) = \left(\int_{S^{n-1}} \|\varphi_\tau(\nabla_uf_\varepsilon)\|_p^{-n}\,du\right)^{-\frac{p}{n}}= \frac{(1+\varepsilon)^n-1}{\varepsilon^{p}} \cdot  \omega_n\cdot  A(n, p).$$   By Definition  \ref{definition1}, for $p\in [1, n)$,   $$\affp(B_n)\leq   \mathcal{H}_{p,\tau}(f_\varepsilon)\Big|_{\varepsilon=1} <  2^n \cdot  \omega_n\cdot  A(n, p)<\infty.$$
We would like to mention that the general $p$-affine capacity can be also defined for $p\in (0, 1)\cup [n, \infty)$ along the same manner in Definition  \ref{definition1}, however in these cases the general $p$-affine capacities are trivial. For instance, if $p\in (0, 1)$,  $$\affp(B_n)\leq \lim_{\varepsilon\rightarrow 0^+} \mathcal{H}_{p,\tau}(f_\varepsilon) =\lim_{\varepsilon\rightarrow 0^+} \frac{(1+\varepsilon)^n-1}{\varepsilon^{p}} \cdot  \omega_n\cdot  A(n, p)=0,$$ and hence, again due to  the proofs of (\ref{monotone-p-affine-capacity}) and Theorem \ref{property},  $\affp(K)=0$ for any compact set $K\subset \Rn$ and for any $\tau\in [-1, 1]$. The case for $p> n$ can be seen intuitively from the above estimate with $\varepsilon\rightarrow \infty$ instead, but more details for  $p\geq n$ will be discussed in  Theorem \ref{comparison}. The precise value of $\affp(B_n)$ will be provided in formulas (\ref{affine-isocapacitary-2-2}) and (\ref{ball-affine-p-1}). $\hfill \Box$

 As $\varphi_0(t)=2^{-1/p} |t|$, one gets  the $p$-affine capacity defined by Xiao in \cite{Xiao2015, Xiao}: 
$$C_{p,0}(K)=\frac{1}{2}\inf_{f\in \mathcal{E}(K)}  \left(\int_{S^{n-1}}\lVert \nabla_u f\lVert_{p}^{-n}\ du\right)^{-\frac{p}{n}}.
$$  As $\varphi_1(\nabla_u f)= \nabla^{+}_u f,$ one has  $$
C_{p,1}(K)= \inf_{f\in \mathcal{E}(K)}\left(\int_{S^{n-1}}\lVert \nabla^{+}_u f\lVert_{p}^{-n}\ du\right)^{-\frac{p}{n}},
$$ which will be called the asymmetric $p$-affine capacity and denoted by $C_{p,+}$ instead of $C_{p, 1}$ for better intuition. Similarly,  as $\varphi_{-1} (\nabla_u f)= \nabla^{-}_u f,$ one can have the following $p$-affine capacity:  
$$C_{p,-}(K)=\inf_{f\in \mathcal{E}(K)} \left(\int_{S^{n-1}}\lVert \nabla^{-}_u f\lVert_{p}^{-n}\ du\right)^{-\frac{p}{n}}.
$$

 The following theorem plays important roles in later context. For a compact set $K\subset \Rn$, let $$\mathcal{F}(K)=\Big\{f: f\in W_0^{1, p}, \ 0\leq f\leq 1 \ \mathrm{in}\ \R^n, \ \mathrm{and}\ f=1 \ \mathrm{in \ a \ neighborhood \ of \ K} \Big\}.$$

\bt \label{definition2} Let $1\leq p<n$ and $K$ be a compact set in $\Rn$. Then 
$$C_{p,\tau}(K)=\inf_{f\in  \mathcal{F}(K)} \mathcal{H}_{p,\tau}(f).$$ Moreover,  the general $p$-affine capacity is upper-semicontinuous:  for any $\varepsilon>0$, there exists an open set $O_{\varepsilon}$ such that for any compact set $F$ with $K\subset F\subset O_\varepsilon$,  $$\affp(F)\leq \affp(K)+\varepsilon.$$
\et  \begin{proof} Our proof is based on the standard technique in \cite{Mazya} and is similar to that in \cite{Xiao,Xiaozhang}. A short proof is included for completeness. Recall that $\affp(K)<\infty$.  Due to  $\mathcal{F}(K)\subset \mathcal{E}(K)$, one has  $$\inf_{f\in  \mathcal{F}(K)} \mathcal{H}_{p,\tau}(f) \geq \affp(K).$$ On the other hand, for any $\varepsilon>0$, let $f_{\varepsilon}\in \mathcal{E}(K)$ satisfy that $$
	C_{p,\tau}(K)+\varepsilon\geq\mathcal{H}_{p,\tau}(f_\varepsilon).
	$$  For  $i=1, 2, \cdots,$ there are functions $\phi_i\in C_c^{\infty}(\R)$,  such that,   for all $t \in \R$, $$0\leq \phi'_i(t)\leq i^{-1}+1, $$    $\phi_i=0$ in  a  neighborhood  of $(-\infty,0]$, and  $\phi_i=1$ in  a  neighborhood of $[1,\infty).$ It follows from the chain rule in \cite[Theorem 4 on p.129]{Evan} and the homogeneity of $\varphi_\tau$ (see (\ref{taufunction}))  that,  for all $i$,   $\phi_i (f_\varepsilon)\in \mathcal{F}(K)$   and  
	\begin{eqnarray*}
		\inf_{f\in  \mathcal{F}(K)} \mathcal{H}_{p,\tau}(f) & \leq & \mathcal{H}_{p,\tau}(\phi_i (f_\varepsilon))\\ &\leq&  (1+i^{-1})^p \cdot \mathcal{H}_{p,\tau}(f_\varepsilon) \\&\leq&  (1+i^{-1})^p \cdot (C_{p,\tau}(K)+\varepsilon).
	\end{eqnarray*}
	Taking $i\rightarrow \infty$ first and then letting $\varepsilon\rightarrow 0$, one gets  $$\inf_{f\in  \mathcal{F}(K)} \mathcal{H}_{p,\tau}(f) \leq  \affp(K)$$ and hence the following desired formula holds: $$\inf_{f\in  \mathcal{F}(K)} \mathcal{H}_{p,\tau}(f)= \affp(K).$$

	Now let us prove the upper-semicontinuity. For any given $\varepsilon>0$, let $g_{\varepsilon}\in \mathcal{F}(K)$ and $O_{\varepsilon}$  be a neighborhood of $K$ such that $g_{\varepsilon}=1$ on $O_{\varepsilon}$ and $$\affp(K)+\varepsilon\geq \mathcal{H}_{p,\tau}(g_{\varepsilon}).$$  On the other hand, for any compact set $F$ such that $K\subset F\subset O_{\varepsilon}$, one has $g_{\varepsilon} \in \mathcal{F}(F)$ and hence  
	$$ \mathcal{H}_{p,\tau}(g_{\varepsilon}) \geq C_{p,\tau}(F),
	$$ by Definition \ref{definition1}. The desired inequality follows from the above two inequalities.  \end{proof}

Our next result regarding the definition of the general $p$-affine capacity for compact sets is to replace $\mathcal{E}(K)$ by the bigger set  $\mathcal{D}(K):$ $$\mathcal{D}(K)=\Big\{ f\in W_0^{1, p}\ \ \mathrm{such\ that}\ \ f\geq 1 \ \mathrm{on}  \ K\Big\}.$$

\bt \label{prop 3-3} Let $1\leq p<n$ and $K$ be a compact set in $\Rn$.  Then  \begin{eqnarray*} 
	\affp(K)&=&\inf_{f\in \mathcal{D}(K)}  \mathcal{H}_{p,\tau}(f).
\end{eqnarray*} \et   
\begin{proof} It follows from (\ref{function-tau}) and \cite[Lemma 1.19]{HKM-1} that, for any $f\in W_0^{1, p}$ and for any $u\in \sphere$, $$\varphi_{\tau}(\nabla_u f_+(x))=  \left\{\begin{array}{ll} \varphi_{\tau}(\nabla_u f(x)), \ \ & \ \ \text{if} \ \  f(x)>0, \\ 0,  \ \ & \ \ \text{if} \ \  f(x)\leq 0.  \end{array} \right.$$  Hence, for any $u\in \sphere$ and all $x\in \Rn$, one has $$\varphi_{\tau}(\nabla _uf_+(x))\leq \varphi_{\tau}(\nabla _uf(x)).$$ This further implies that $\mathcal{H}_{p,\tau}(f_+)\leq \mathcal{H}_{p,\tau}(f)$  for any $f\in W_0^{1, p}$.  Let $\{f_k\}_{k\geq 1}\subset  \mathcal{D}(K)$ be  such that  
	$$\lim_{k\rightarrow \infty}\mathcal{H}_{p,\tau}(f_k)=
	\inf_{f\in \mathcal{D}(K)}  \mathcal{H}_{p,\tau}(f).$$ Then $\{f_{k,+}\}_{k\geq 1}$ is a sequence in $\mathcal{E}(K)$. Definition \ref{definition1} yields \begin{eqnarray*}
		\lim_{k\rightarrow \infty}\mathcal{H}_{p,\tau}(f_k) &\geq & \limsup_{k\rightarrow \infty}\mathcal{H}_{p,\tau}(f_{k,+})\geq \inf_{f\in \mathcal{E}(K)}\mathcal{H}_{p,\tau}(f)= \affp(K).\end{eqnarray*} This concludes that $$\inf_{f\in \mathcal{D}(K)}  \mathcal{H}_{p,\tau}(f)\geq \affp(K).$$ On the other hand, as $\mathcal{E}(K)\subset \mathcal{D}(K)$, the following inequality holds trivially:  $$
	\inf_{f\in \mathcal{D}(K)}  \mathcal{H}_{p,\tau}(f)\leq \affp(K).$$ Combining the above two inequalities, one has  $\affp(K)=\inf_{f\in \mathcal{D}(K)}  \mathcal{H}_{p,\tau}(f).$  \end{proof}

The following result asserts that  $f \in W_0^{1, p}$ in Definition \ref{definition1}, Theorems \ref{definition2} and \ref{prop 3-3}  could be replaced by $f \in C_c^{\infty}$. The smoothness of functions is convenient in establishing many properties for the general $p$-affine capacity. 

\bt\label{definition-smooth-1} Let $p\in [1, n)$ and $K$ be a compact set in $\Rn$. For any $\tau\in [-1, 1]$, one has   \begin{eqnarray} 
\affp(K)=\inf_{ f\in C_c^{\infty}\cap \mathcal{D}(K)} \mathcal{H}_{p,\tau}(f) = \inf_{f\in C_c^{\infty}\cap \mathcal{E}(K) } \mathcal{H}_{p,\tau}(f)   =\inf_{f\in C_c^{\infty}\cap \mathcal{F}(K)}  \mathcal{H}_{p,\tau}(f).\label{replaced-by-smooth-functions-11}
\end{eqnarray} \et  \begin{proof} Let $p\in [1, n)$.  Let $f\in \mathcal{F}(K)$,  i.e., $f\in W_0^{1, p}$ such that $0\leq f\leq 1$  in $\R^n$ and $f=1$ in $U$, a  neighborhood of $K$. As $W_0^{1, p}$ is the closure of $C_{c}^{\infty}$ under $\|\cdot\|_{1, p}$, there is a sequence $\{f_k\}_{k=1}^{\infty} \subset C_{c}^{\infty}$  such that $f_k\rightarrow f$ in $W_0^{1, p}$, i.e., $$ \| f_k-f\|_p+\| \nabla f_k-\nabla f\|_p\rightarrow 0.$$ Without loss of generality, we can assume that $f_k\in C_c^{\infty}\cap \mathcal{D}(K)$ for all $k.$ To see this, from the regularization technique (see, e.g., \cite{HKM-1}), one can choose a cut off function $\kappa\in C^{\infty}$, such that, $0\leq \kappa \leq 1$ on $\Rn$,  $\kappa=1$ on  $\mathbb{R}^n\setminus U,$ and $\kappa=0$ in a neighborhood (contained in $U$) of $K$. Let $$g_k=1-(1-f_k)\kappa.$$ Clearly,  $g_k\in C_c^{\infty}$, such that, $g_k=1$ in  a neighborhood (contained in $U$) of $K$ and $g_k=f_k$ on $\Rn\setminus U$. This implies $g_k\in C_c^{\infty}\cap \mathcal{D}(K)$ for all $k$.  Moreover, $\|g_k-f\|_{1, p}\rightarrow 0$ and hence $$\|g_k-f\|_p\rightarrow 0 \ \ \ \mathrm{and}\ \ \ \|\nabla g_k-\nabla f\|_p\rightarrow 0. $$

 Let $f_k\in C_c^{\infty}\cap \mathcal{D}(K)$  be such that $f_k\rightarrow f$ in $W_0^{1, p}$. It can be checked that, for any $u\in \sphere$,  $$|\nabla^+_uf_k-\nabla^+_uf|\leq |\nabla f_k-\nabla f| \ \ \ \mathrm{and}\ \ \ \ |\nabla^-_uf_k-\nabla^-_uf|\leq |\nabla f_k-\nabla f|. $$ This together with (\ref{function-tau}) yield, for any $\tau\in [-1, 1]$ and for all $k\geq 1$,  \begin{eqnarray*} \Big|\varphi_\tau(\nabla_uf_k)- \varphi_\tau(\nabla_uf)\Big|\!\!&=&\!\!\Big| \Big(\frac{1+\tau}{2}\Big)^{1/p}  \big[\nabla_u^+f_k-\nabla_u^+f\big] +\Big(\frac{1-\tau}{2}\Big)^{1/p}  \big[\nabla_u^-f_k-\nabla_u^-f\big] \Big|\\\!\!&\leq&\!\! \Big(\frac{1+\tau}{2}\Big)^{1/p}  \Big| \nabla_u^+f_k-\nabla_u^+f\Big| +\Big(\frac{1-\tau}{2}\Big)^{1/p}  \Big|\nabla_u^-f_k-\nabla_u^-f\Big|\\ \!\!&\leq &\!\!   C(p, \tau) \cdot \big|\nabla f_k-\nabla f\big|, \end{eqnarray*} where we have let $C(p, \tau)$ be the constant $$C(p, \tau) = \Big(\frac{1+\tau}{2}\Big)^{1/p} +\Big(\frac{1-\tau}{2}\Big)^{1/p}.$$ It follows from the triangle inequality that, for any $u\in \sphere$, for any $\tau\in [-1, 1]$  and for any $p\in [1, n)$,  \begin{eqnarray*}  \Big|\|\varphi_\tau(\nabla_uf_k)\|_p-\| \varphi_\tau(\nabla_uf)\|_p \Big|&\leq&   \|\varphi_\tau(\nabla_uf_k)- \varphi_\tau(\nabla_uf)\|_p\\ &\leq& C(p, \tau) \cdot \|\nabla f_k-\nabla f\|_p. \end{eqnarray*}     Consequently,  for any $u\in \sphere$, for any $\tau\in [-1, 1]$  and for any $p\in [1, n)$,  one has   
	\begin{equation*}  
	\lim_{k\rightarrow\infty} \| \varphi_\tau(\nabla_uf_k)\|_p=\| \varphi_\tau(\nabla_uf)\|_p.
	\end{equation*}  By Fatou's lemma, one has 
	\begin{eqnarray} \mathcal{H}_{p,\tau}(f)&=&  \left(\int_{S^{n-1}} \|\varphi_\tau(\nabla_uf)\|_p^{-n}\,du\right)^{-\frac{p}{n}} \nonumber\\
	&=&  \left(\int_{S^{n-1}}\lim_{k\rightarrow\infty} \|\varphi_\tau(\nabla_uf_k)\|_p^{-n}\,du\right)^{-\frac{p}{n}} \nonumber\\ &\geq& \left(\liminf_{k\rightarrow\infty}\int_{S^{n-1}} \|\varphi_\tau(\nabla_uf_k)\|_p^{-n}\,du\right)^{-\frac{p}{n}} \nonumber \\&=&  \limsup_{k\rightarrow\infty}\left(\int_{S^{n-1}} \|\varphi_\tau(\nabla_uf_k)\|_p^{-n}\,du\right)^{-\frac{p}{n}} \nonumber \\&=& \limsup_{k\rightarrow\infty} \mathcal{H}_{p,\tau}(f)\nonumber \\ &\geq & \inf_{g\in C_c^{\infty}\cap \mathcal{D}(K)} \mathcal{H}_{p,\tau}(g).\label{continuous-approximation-1}
	\end{eqnarray} It follows from Theorem \ref{definition2}  that, by taking the infimum over  $f\in \mathcal{F}(K)$,   \begin{eqnarray*}
	\affp(K)  \geq \inf_{C_c^{\infty}\cap \mathcal{D}(K)} \mathcal{H}_{p,\tau}(f).
\end{eqnarray*} It is easily checked that, due to $C_c^{\infty}\subset W_0^{1, p}$,  \begin{eqnarray*} 
\affp(K)\leq \inf_{C_c^{\infty}\cap \mathcal{D}(K)} \mathcal{H}_{p,\tau}(f),
\end{eqnarray*} and hence equality holds, as desired.  

The desired formula (\ref{replaced-by-smooth-functions-11}) follows, due to $\mathcal{F}(K)\subset \mathcal{E}(K)\subset \mathcal{D}(K)$, once the following inequality is proved: $$
\inf_{f\in C_c^{\infty}\cap \mathcal{F}(K)}  \mathcal{H}_{p,\tau}(f)\leq \inf_{f\in C_c^{\infty}\cap \mathcal{D}(K)}  \mathcal{H}_{p,\tau}(f)=\affp(K).$$ This inequality follows along the same lines as the proof of Theorem \ref{definition2}. In fact, for any $\varepsilon>0$, let $f_{\varepsilon}\in \mathcal{D}(K)\cap C_c^{\infty}$ satisfy that $$
C_{p,\tau}(K)+\varepsilon\geq\mathcal{H}_{p,\tau}(f_\varepsilon).
$$  Let $\phi_i\in C_c^{\infty}(\R)$ be as in Theorem \ref{definition2}. Then, $\phi_i (f_\varepsilon)\in \mathcal{F}(K)\cap C_c^{\infty}$ and  
\begin{eqnarray*}
	\inf_{f\in  \mathcal{F}(K)\cap C_c^{\infty}} \mathcal{H}_{p,\tau}(f)  \leq    (1+i^{-1})^p \cdot (C_{p,\tau}(K)+\varepsilon).
\end{eqnarray*}
Taking $i\rightarrow \infty$ first and then letting $\varepsilon\rightarrow 0$, one gets $$\inf_{f\in   \mathcal{F}(K)\cap C_c^{\infty}} \mathcal{H}_{p,\tau}(f) \leq  \affp(K)$$ as desired.  \end{proof} 

It follows from  (\ref{taufunction}) and  $\nabla_y f=y\cdot \nabla f$ that,  for all $\lambda>0$ and $y\in \Rn\setminus\{o\}$,   $$\|\varphi_\tau(\nabla_{\lambda y}f)\|_p=\lambda \|\varphi_\tau(\nabla_{y}f)\|_p.$$ Moreover,  for $p\in [1, n)$ and for any $y_1, y_2\in  \Rn\setminus\{o\}$, by the Minkowski's inequality, one has 
\begin{eqnarray*}
	\|\varphi_\tau(\nabla_{y_1+y_2}f)\|_p &\leq&  \|\varphi_\tau(\nabla_{y_1} f)+\varphi_\tau(\nabla_{y_2} f)\|_p\\ 
	&\leq&  \|\varphi_\tau(\nabla_{y_1} f)\|_p+\|\varphi_\tau(\nabla_{y_2} f)\|_p.
\end{eqnarray*} Hence, $\|\varphi_\tau(\nabla_{y} f)\|_p: \Rn\setminus\{o\}\rightarrow [0, \infty)$, as a function of $y\in \Rn\setminus\{o\}$, is sublinear.   According to  the proof of  \cite[Lemma 3.1]{VHN} (or \cite[Lemma 2]{HSasymmetric}),  if $f\in \mathcal{F}(K)$, then $\|\varphi_{\tau}(\nabla_u f)\|_p>0$ and $\|\varphi_{\tau}(\nabla_y f)\|_p$  is the support function of a convex body in $\cK_0$. Let $L_{f, \tau}$ be the convex body.  An application of (\ref{supportradial}) and (\ref{volumeformula})  yields  (see also \cite[(3.2)]{VHN}) \begin{eqnarray*}  \mathcal{H}_{p,\tau}(f) &=& \left(\int_{S^{n-1}} \|\varphi_\tau(\nabla_uf)\|_p^{-n}\,du\right)^{-\frac{p}{n}}\\ &=& \left(\int_{S^{n-1}}\big[ h_{L_{f, \tau}}(u)\big]^{-n}\,du\right)^{-\frac{p}{n}} \\ &=& \left(\frac{1}{nV(B_n)}\int_{S^{n-1}}\big[\rho_{L^\circ _{f, \tau}}(u)\big]^{n}\,d\sigma(u)\right)^{-\frac{p}{n}} \\&=& \bigg(\frac{V(L_{f,\tau}^\circ)}{V(B_n)}\bigg)^{-\frac{p}{n}}.\end{eqnarray*}  Taking the infimum over  $f\in \mathcal{F}(K)$,  Theorem \ref{definition2} implies that for any compact set $K\subset \Rn$, for any $\tau\in [-1, 1]$  and for any $p\in [1, n)$,  \begin{eqnarray*} 
	C_{p,\tau}(K)=\inf_{f\in  \mathcal{F}(K)} \mathcal{H}_{p,\tau}(f)= \inf _{ f\in \mathcal{F}(K)}\bigg(\frac{V(L_{f,\tau}^\circ)}{V(B_n)}\bigg)^{-\frac{p}{n}}.
\end{eqnarray*} This provides a connection of the general $p$-affine capacity with the volume of convex bodies.

 The general $p$-affine capacity of a general bounded measurable set $E\subset \Rn$ can be defined as well. In fact, for $O\subset \Rn$ a bounded open set, $$\affp(O)= \sup \Big\{\affp(K): \ \ K\subset O\ \mathrm{and} \ K\ \mathrm{ is \ compact} \Big\}.$$ Then the general $p$-affine capacity of a bounded measurable set $E\subset \Rn$ is formulated by $$\affp(E)=\inf \Big\{\affp(O):  \ \ E\subset O\ \mathrm{and} \ O\ \mathrm{ is \ open} \Big\}.$$ In later context of this article, we only concentrate on the general $p$-affine capacity for compact sets. We would like to mention that many properties proved in Section \ref{section:4}, such as, monotonicity,  affine invariance and homogeneity etc, for compact sets could work for general sets too.

\section{Properties of the general $p$-affine capacity} \label{section:4} This section aims to establish basic properties for the general $p$-affine capacity, such as, monotonicity,  affine invariance, translation invariance, homogeneity and the continuity from above. 

  The following result provides the properties of $\affp(\cdot)$ as a function of $\tau\in [-1, 1]$.  

\bc\label{properties-even-affine-1} Let $p\in [1, n)$ and $K$ be a compact set in $\Rn$. The following properties hold. 
\begin{itemize} \item [i)] For any $\tau\in [-1, 1],$ one has $$\affp(K)=C_{p, -\tau}(K).$$ 
	\item [ii)] For any $\lambda\in [0, 1]$ and for any $\tau, \gamma\in [-1, 1]$,  one has $$C_{p, \lambda \tau+(1-\lambda)\gamma}(K)\geq \lambda \cdot C_{p,  \tau}(K)+(1-\lambda)\cdot C_{p, \gamma}(K).$$ 
\end{itemize}
\ec  \begin{proof} i)    Let $v=-u$. Then for any $x\in \Rn$, one has  \begin{equation*} \nabla^+_u f(x)=\nabla_v^-f(x)\ \ \ \ \mathrm{and} \ \ \ \nabla^-_u f(x)=\nabla_v^+f(x).\end{equation*}  This leads to, as $\,du=\,dv$, for any $f\in \mathcal{E}(K)$,   \begin{eqnarray}  \mathcal{H}_{p,\tau}(f) &=& \left(\int_{S^{n-1}} \|\varphi_\tau(\nabla_uf)\|_p^{-n}\,du\right)^{-\frac{p}{n}} \nonumber  \\ &=&  \left(\int_{S^{n-1}} \left(\int_{\R^n} \Big[\Big(\frac{1+\tau}{2}\Big) (\nabla^+_uf(x))^p +\Big(\frac{1-\tau}{2}\Big) (\nabla^-_uf(x))^p\Big]\ dx\right)^{-\frac{n}{p}}\,du\right)^{-\frac{p}{n}} \nonumber  \\  &=& \left(\int_{S^{n-1}} \left(\int_{\R^n} \Big[\Big(\frac{1+\tau}{2}\Big) (\nabla^-_vf(x))^p +\Big(\frac{1-\tau}{2}\Big) (\nabla^+_vf(x))^p\Big]\ dx\right)^{-\frac{n}{p}}\,dv\right)^{-\frac{p}{n}} \nonumber   \\ &=&   \mathcal{H}_{p,-\tau}(f).\nonumber
	\end{eqnarray} 
	It follows from  Definition \ref{definition1} that, for any $\tau\in [-1, 1]$, for any $p\in [1, n)$ and for any compact set $K\subset \Rn$,  $$\affp(K)=C_{p, -\tau}(K).$$   
	
	ii) For any $\lambda\in [0, 1]$ and for any $\tau, \gamma\in [-1, 1],$ it follows from (\ref{function-tau-p}) that, for any $t\in \R$, 
	\begin{equation*} 
	\big[\varphi_{\lambda \tau+(1-\lambda)\gamma}(t)\big]^p =\lambda \big[\varphi_{\tau}(t)\big]^p +(1-\lambda)\big[\varphi_{\gamma}(t)\big]^p,
	\end{equation*} which implies  \begin{eqnarray*}   \int_{\Rn} \big[\varphi_{\lambda \tau+(1-\lambda)\gamma}(\nabla_uf(x))\big]^p\,dx= \lambda \int_{\Rn} \big[\varphi_{\tau}(\nabla_uf(x))\big]^p\,dx+ (1-\lambda) \int_{\Rn} \big[\varphi_{\gamma}(\nabla_uf(x))\big]^p\,dx.
	\end{eqnarray*} According to  the proof of  \cite[Lemma 3.1]{VHN} (or \cite[Lemma 2]{HSasymmetric}), $\|\varphi_{\tau}(\nabla_u f)\|_p>0$ if $f\in \mathcal{F}(K)$. The reverse Minkowski inequality yields that  for any $\lambda\in [0, 1]$ and for any $\tau, \gamma\in [-1, 1],$
	\begin{eqnarray*}  \left(\int_{S^{n-1}}\! \! \|\varphi_{\lambda \tau+(1-\lambda)\gamma}(\nabla_uf)\|_p^{-n}\,du\!\right)^{-\frac{p}{n}} \! \geq \lambda \left(\int_{S^{n-1}}\!\! \|\varphi_{\tau}(\nabla_uf)\|_p^{-n}\,du\!\right)^{-\frac{p}{n}} \!\!+  (1-\lambda) \left(\int_{S^{n-1}}\!\! \|\varphi_{\gamma}(\nabla_uf)\|_p^{-n}\,du\!\right)^{-\frac{p}{n}}. \end{eqnarray*} 
	Taking the infimum over $f\in \mathcal{F}(K)$, by Theorem \ref{definition2}, $$C_{p, \lambda \tau+(1-\lambda)\gamma}(K)\geq \lambda \cdot C_{p,  \tau}(K)+(1-\lambda)\cdot C_{p, \gamma}(K) $$ holds  for any $\lambda\in [0, 1]$ and for any $\tau, \gamma\in [-1, 1].$  	 \end{proof} 

From Corollary \ref{properties-even-affine-1}, one sees that, for any $p\in [1, n)$ and for any compact set $K\subset \Rn$,  $\affp(K)\leq C_{p, \gamma}(K)$ holds  if $-1\leq \tau<\gamma \leq 0,$ and   $C_{p, \gamma}(K)\leq \affp(K)$ holds if  $0\leq \tau< \gamma\leq 1.$ In particular, for any $\tau\in [-1, 1]$, one has $$C_{p, +}(K)=C_{p, -}(K)\leq \affp(K)=C_{p,-\tau}(K)\leq C_{p, 0}(K). $$

  Given two compact sets $K\subset L$, one sees $\mathcal{E}(L)\subset \mathcal{E}(K)$ and hence the general $p$-affine capacity is monotone by Definition \ref{definition1}, namely, \begin{equation}\label{monotone-p-affine-capacity} \affp (K)\leq \affp(L).\end{equation} The general $p$-affine capacity is also translation invariant.  To see this, let $a\in \Rn$ and consider the function  $g(x)=f(x+a)$ for any $x\in \Rn$. It is easily checked that  $f\in \mathcal{E}(K+a)$ if and only if $g\in \mathcal{E}(K)$. Moreover,  $\nabla g(x)=\nabla f(x+a)$, and thus $\mathcal{H}_{p,\tau}(g)=\mathcal{H}_{p,\tau}(f).$ Taking the infimum over $g\in  \mathcal{E}(K)$ from both sides, by Definition \ref{definition1}, for any $a\in \R^n$ and for any compact set $K\subset \Rn$, 
$$\affp (K+a)=\affp(K).$$  An interesting (and common for many capacities) fact for the general $p$-affine capacity is that $$\affp(K)=\affp(\partial K)$$ for any compact set $K\subset \Rn$. To see this, let  $\varepsilon>0$ be given. There exists $f_{\varepsilon}\in \mathcal{E}(\partial K)$ such that
$$
C_{p,\tau}(\partial K)+\varepsilon\geq \mathcal{H}_{p,\tau}(f_{\varepsilon}).
$$
Let $g=\max\{f_{\varepsilon}, 1\}$ on $K$ and $g=f_{\varepsilon}$ on $\Rn\setminus K$.  It can be checked, along the manner same as the proof of Theorem \ref{prop 3-3}, that $g\in \mathcal{E}(K)$ and 
$$
\int_{\Rn}\big[\varphi_\tau(\nabla_u g)\big]^p\ dx \leq \int_{\Rn}\big[\varphi_\tau(\nabla_u f_{\varepsilon})\big]^p\ dx.
$$ Consequently, due to Definition \ref{definition1},  
\begin{eqnarray*}
	C_{p,\tau}(K) \leq \mathcal{H}_{p,\tau}(g) 
	\leq \mathcal{H}_{p,\tau}(f_{\varepsilon}) <C_{p,\tau}(\partial K)+\varepsilon.
\end{eqnarray*}
Letting $\varepsilon\rightarrow 0$, one gets $$C_{p,\tau}(K)\leq C_{p,\tau}(\partial K).$$  The monotonicity of the general $p$-affine capacity yields that $$C_{p,\tau}(\partial K)\leq C_{p,\tau}(K)$$ and hence $C_{p,\tau}(\partial K) = C_{p,\tau}(K)$ holds for all compact set $K\subset \Rn$.   

  Let $GL(n)$ be the group of all invertible linear transforms defined on $\Rn$. For $T\in GL(n)$, denote by $T^t$ and $\det(T)$ the transpose of $T$ and the determinant of $T$, respectively. The affine invariance of the general $p$-affine capacity is stated in the following theorem.   \bt \label{property}
The general $p$-affine capacity has the affine invariance and homogeneity: for any $T\in GL(n)$ and for any compact set $K\subset \Rn$, 
$$\affp (TK)=|\det(T)|^{\frac{n-p}{n}} \affp(K).  $$  In particular,  the general $p$-affine capacity is affine invariant: for any $T\in GL(n)$ with $|\det(T)|=1$, 
$$\affp (TK)= \affp(K).$$ Moreover,  the general $p$-affine capacity has positive homogeneity of degree $n-p$, i.e., $$\affp(\lambda K)=\lambda^{n-p} \affp(K)$$ for all $\lambda>0$, where $\lambda K=\{\lambda x: x\in K\}$.  
\et 

\begin{proof}   For  $T\in GL(n)$ and $f\in \mathcal{E}(TK)$, one has  $g=f\circ T\in \mathcal{E}(K)$. For simplicity, assume that $|\det(T)|=1$. Thus, by $x=Ty$,
	\begin{eqnarray*} \int_{\R^n}\big[\varphi_\tau(\nabla_u g(y))\big]^p \,dy
		=\int_{\R^n}\big[\varphi_\tau(\nabla_u (f\circ T)(y))\big]^p \,dy
		=\int_{\R^n} \big[\varphi_\tau(\nabla_{Tu} (f(x))\big]^p \,dx,
	\end{eqnarray*} where the second equality follows from the chain rule  
	$$\nabla g(y)=
	\nabla(f\circ T)(y)=T^t\nabla f(Ty).
	$$ By letting $v=Tu/|Tu|$, it follows from  (\ref{taufunction}) that 
	\begin{eqnarray*}
		\int_{S^{n-1}} \left(\int_{\R^n}\big[\varphi_\tau(\nabla_u g(y))\big]^p \,dy\right)^{-\frac{n}{p}}\, du &=& \int_{S^{n-1}} \left(\int_{\R^n}\big[\varphi_\tau(\nabla_{Tu} (f(x))\big]^p \,dx\right)^{-\frac{n}{p}}\,du\\
		&=& \int_{S^{n-1}} \left(\int_{\R^n}\big[\varphi_\tau(\nabla_{v} (f(x))\big]^p  \,dx\right)^{-\frac{n}{p}}|Tu|^{-n}\,du\\
		&=&\int_{S^{n-1}} \left(\int_{\R^n} \big[\varphi_\tau(\nabla_{v} (f(x))\big]^p  \,dx\right)^{-\frac{n}{p}}\,dv.
	\end{eqnarray*} Consequently, $\mathcal{H}_{p,\tau}(g)=\mathcal{H}_{p,\tau}(f)$. 
	Taking the infimum over $f\in \mathcal{E}(TK)$ from both sides, which is equivalent to taking the infimum over $g\in \mathcal{E}(K)$ from the left hand side, one gets the affine invariance: for all $T\in GL(n)$ with $|\det (T)|=1$, then   $$C_{p,\tau}(TK)=C_{p,\tau}(K).$$ 
	
	For the homogeneity,  let $\lambda>0$ be given. For any $f\in \mathcal{E}(\lambda K)$, one sees $g_{\lambda}\geq \mathbf{1}_K$ where  $g_\lambda(x)=f(\lambda x)$ for all $x\in \Rn$. 
	It is easily checked, by letting $y=\lambda x$, that $$
	\int_{\R^n} \big[\varphi_\tau(\nabla_u g_\lambda (x))\big]^p\, dx=\lambda^{p-n}\int_{\R^n} \big[\varphi_\tau(\nabla_uf(y))\big]^p\, dy,
	$$ which further implies that $ \mathcal{H}_{p,\tau}(f)=
	\lambda^{n-p} \mathcal{H}_{p,\tau}(g_\lambda).$  The desired formula $\affp(\lambda K)=\lambda^{n-p} \affp(K)$  follows immediately by Definition \ref{definition1} and by taking the infimum over $f\in \mathcal{E}(\lambda K)$.

  Finally, we consider $T\in GL(n)$ be an invertible linear transform. Then $$\widetilde{T}=|\det(T)|^{-1/n} T$$ has $|\det (\widetilde{T})|=1$.   Hence, the affine invariance and the homogeneity yield that, for all $T\in GL(n)$, 
	$$\affp (TK) =\affp(|\det(T)|^{1/n}\widetilde{T}K  )=|\det(T)|^{\frac{n-p}{n}} \affp(\widetilde{T}K)=|\det(T)|^{\frac{n-p}{n}} \affp(K).$$ This concludes the proof. \end{proof}

The continuity from above for the general $p$-affine capacity is stated in the following theorem. 

\bt The general $p$-affine capacity is continuous from above:  if $\{K_i\}_{i=1}^{\infty}$ is a decreasing sequence of compact sets, then \begin{equation} \label{continuous from above} \affp(\cap_{i=1}^{\infty} K_i)=\lim_{i\rightarrow \infty} \affp(K_i).\end{equation} 
\et  
\begin{proof} Recall that the general $p$-affine capacity of the compact set $K_1$ is finite. It follows from the monotonicity that,  for all $i$, $$\affp(K_{i+1})\leq \affp(K_i)\leq \affp(K_1)<\infty,$$ and hence $\lim_{i\rightarrow \infty}\affp(K_i)$ exists and is finite. Moreover, the monotonicity of the general $p$-affine capacity also yields  $$\affp(\cap _{i=1}^{\infty} K_i) \leq \lim_{i\rightarrow \infty}\affp(K_i).$$ 
	
	The desired formula (\ref{continuous from above}) follows if we prove the following inequality:  $$\affp(\cap _{i=1}^{\infty} K_i) \geq \lim_{i\rightarrow \infty}\affp(K_i).$$  First of all, the set $\cap_{i=1}^{\infty}  K_i$ is clearly compact. By Definition \ref{definition1} and Theorem \ref{definition-smooth-1},  for any $\varepsilon>0$, one can find a smooth function $f_{\varepsilon} \in \mathcal{E}(\cap_{i=1}^{\infty} K_i )$, such that,
	$f_{\varepsilon}\geq  \mathbf{1}_{\cap_{i=1}^{\infty} K_i }$
	and $$
	\affp\big(\!\cap_{i=1}^{\infty}  K_i\big)+\varepsilon\geq \mathcal{H}_{p,\tau}(f_{\varepsilon}).$$
	Let $K_\varepsilon=\{x\in \R^n: f_{\varepsilon}(x)\geq 1-\varepsilon \}.$  Then, $\frac{f_{\varepsilon}}{1-\varepsilon}\in \mathcal{E}(K_{\varepsilon})$ and  $K_i\subset K_\varepsilon$  for $i$ big enough. Together with (\ref{taufunction}), Definition  \ref{definition1} and the monotonicity of the general $p$-affine capacity, one has  
	\begin{eqnarray*}
		\lim_{i\rightarrow\infty}\affp(K_i) \leq \affp(K_\varepsilon) \leq (1-\varepsilon)^{-p} \mathcal{H}_{p,\tau}(f_\varepsilon)\leq \frac{\affp(\cap_{i=1}^{\infty} K_i)+\varepsilon}{(1-\varepsilon)^{p}}.
	\end{eqnarray*}  Taking $\varepsilon\rightarrow 0$, one gets the desired inequality $$\lim_{i\rightarrow\infty}\affp(K_i) \leq  \affp(\cap_{i=1}^{\infty} K_i)$$ and this concludes the proof. \end{proof} 

Note that one cannot expect to have the subadditivity for the general $p$-affine capacity, even for $\tau=0$;  see \cite{Xiaozhang} for the details. It is not clear whether the general $p$-affine capacity has the continuity from below.

\section{Sharp geometric inequalities for the general $p$-affine capacity}\label{section:5}
This section aims to establish several sharp geometric inequalities for the general $p$-affine capacity. In particular, the general $p$-affine capacity is compared with the $p$-variational capacity, the general $p$-integral affine surface areas and the volume. 

\subsection{Comparison with the $p$-variational capacity} 
This subsection aims to compare the general $p$-affine capacity and the $p$-variational capacity. For $p\in [1, n)$ and a compact set $K\subset \Rn$, the $p$-variational capacity of $K$, denoted by $C_p(K)$, is formulated by   \begin{eqnarray*} C_p(K) &=&\inf_{f\in \mathcal{D}(K)}  \int_{\R^n} |\nabla f|^p\,dx= \inf_{f\in \mathcal{D}(K) \cap C_c^{\infty}}  \int_{\R^n} |\nabla f|^p\,dx.\end{eqnarray*}  Of course, the set $\mathcal{D}(K)$ in the above definition for the $p$-variational capacity could be replaced by $\mathcal{E}(K)$ and $\mathcal{F}(K)$ (see e.g., \cite{Evan, Mazya}).  The $p$-variational capacity is fundamental in many areas, such as, analysis, geometry and physics. It has many properties similar to those for the general $p$-affine capacity, such as, homogeneity, monotonicity; however the $p$-variational capacity does not have the affine invariance. 

The comparison between the general $p$-affine capacity and the $p$-variational capacity is stated in the following theorem. The case $\tau=0$ was discussed in \cite[Remark 2.7]{Xiao} and  \cite[Theorem 1.5']{Xiao2015}. Let $A(n, p)$ be the constant given in (\ref{computation}). 
\bt \label{comparison} Let $p\in [1, n)$ and $K\subset \Rn$ be a compact set. For any $\tau\in [-1, 1]$, one has 
$$
C_{p,\tau}(K)\leq A(n, p) \cdot  C_{p}(K). 
$$  
\et

\begin{proof}  According to  the proof of  \cite[Lemma 3.1]{VHN} (or \cite[Lemma 2]{HSasymmetric}),   $\|\varphi_{\tau}(\nabla_u f)\|_p>0$  for any $f\in \mathcal{F}(K)\cap C_c^{\infty}$,  for any $\tau\in [-1, 1]$ and for any $u\in \sphere$. By Jensen's inequality, Fubini's theorem, (\ref{taufunction}) and (\ref{computation}), one has, for any $f\in \mathcal{F}(K)\cap C_c^{\infty}$, 
	\begin{eqnarray}\mathcal{H}_{p,\tau}(f)&=&
	\bigg(\int_{S^{n-1}} \bigg(\int_{\R^n}\big[\varphi_\tau(\nabla_uf)\big]^p\, dx\bigg)^{-\frac{n}{p}}\, du\bigg)^{-\frac{p}{n}}\nonumber \\ &\leq& \int_{S^{n-1}} \bigg(\int_{\R^n}\big[\varphi_\tau(\nabla_uf)\big]^p\, dx \bigg)\, du\nonumber \\   
	&=&\int_{\Rn} \bigg( \int_{S^{n-1}} \big[\varphi_\tau(\nabla_uf)\big]^p\, du\bigg) \, dx\nonumber\\ 
	&= &\bigg( \int_{S^{n-1}}\big[\varphi_\tau(u\cdot v)\big]^p \,du \bigg)\cdot \bigg(\int_{\R^n}|\nabla f|^p \ dx\bigg) \nonumber\\
	&=&A(n, p)\cdot \int_{\R^n}|\nabla f|^p \ dx, \nonumber
	\end{eqnarray} where  $v\in \sphere$ (depending on $x\in \Rn$) is given by
	$$
	v=\frac{\nabla f(x)}{|\nabla f(x)|} \ \ \  \mathrm{on} \ \ \ \{x\in \Rn: \nabla f\neq 0\}. 
	$$ Taking the infimum over $f\in \mathcal{F}(K)\cap C_c^{\infty}$, one has, by Theorem \ref{definition-smooth-1} and the definition of the $p$-variational capacity,  
	\begin{eqnarray*} \affp(K)&=&\inf_{f\in \mathcal{F}(K)\cap C_c^{\infty}} \mathcal{H}_{p,\tau}(f) \\ &\leq&  A(n, p)\cdot \inf_{f\in \mathcal{F}(K)\cap C_c^{\infty}} \int_{\R^n}|\nabla f|^p \ dx\\ &=&A(n, p)\cdot  C_p(K)
	\end{eqnarray*} holds for any $\tau\in [-1, 1]$, for any $p\in [1, n)$ and for any compact set $K\subset \Rn$.  \end{proof}

  It is well known  (see e.g., \cite[(2.2.13) and (2.2.14)]{Mazya}) that \begin{equation} \label{capacity-ball-1} C_p(B_n)=n\omega_n\cdot \left(\frac{n-p}{p-1}\right)^{p-1} \end{equation}  for $p\in (1, n)$,  $C_p(B_n)=0$ for $p\geq n$, and $  C_1(B_n)=\lim_{p\rightarrow 1^+} C_p(B_n)=n\omega_n.$ Hence,  for any $\tau\in [-1, 1]$,  \begin{eqnarray} \affp(B_n)\leq  A(n, p)C_p(B_n)= A(n, p)\cdot  n\omega_n\cdot \left(\frac{n-p}{p-1}\right)^{p-1}\label{inequality-affine-capacity-ball-1}\end{eqnarray} holds for any $p\in (1, n)$,  and \begin{eqnarray} C_{1, \tau}(B_n)\leq  A(n, 1)C_1(B_n)= A(n, 1)\cdot  n\omega_n.\label{inequality-affine-capacity-ball-2-2}\end{eqnarray}
Following along the same lines as the proof of Theorem \ref{comparison}, one has, for any $\tau\in [-1, 1]$ and for any $p\geq n$,  $$0\leq C_{p,\tau}(B_n)\leq A(n, p)C_p(B_n)=0.$$ Again due to  the proofs of (\ref{monotone-p-affine-capacity}) and Theorem \ref{property},  $C_{p, \tau}(K)=0$ for any $\tau\in [-1, 1]$,  for any $p\geq n$ and for any compact set $K\subset \Rn$. 
\subsection{Affine isocapacitary inequalities} 
This subsection dedicates to establish the affine isocapacitary inequality which compares the general $p$-affine capacity with the volume. An ellipsoid is a convex body of form $TB_n+x_0$ for some $T\in GL(n)$ and $x_0\in \Rn$.

\bt \label{newisocapacity} Let $p\in [1, n)$. For any $\tau\in [-1, 1]$ and for any compact set $K\subset \Rn$, one has $$ \bigg(\frac{\affp(K)}{\affp(B_n)}\bigg)^{\frac{1}{n-p}}\geq \bigg(\frac{V(K)}{V(B_n)}\bigg)^{\frac{1}{n}} $$ 
with equality if $K$ is an ellipsoid.
\et

\begin{proof} Let $p\in (1, n)$, $\tau\in [-1, 1]$ and $K\subset \Rn$ be a compact set. It follows from \cite[inequality (5.8)] {HSasymmetric} that for $f\in C_c^{\infty} \cap \mathcal{F}(K)$,  $\|f\|_{\infty}=1$ and  $$
	\left(\int_{S^{n-1}}\| \nabla^{+}_u f\|_{p}^{-n}\ du\right)^{-\frac{p}{n}}\geq n^p\omega_n^{\frac{p}{n}}A(n,p) \int_{0}^{1}\frac{V([f]_t)^{\frac{np-p}{n}}}{[-V([f]_t)']^{p-1}} \ dt,
	$$ where $V([f]_t)'$ is the derivative of $V([f]_t)$ with respect to $t$. Recall that for any real number $t>0$ and for any $f\in C_{c}^{\infty}$, 
	$$
	[f]_t=\{x\in \R^n: |f(x)|\geq t\}.$$  Note that $V(K)\leq V([f]_1)\leq V([f]_0)$. Together with 
	Jensen's inequality, one has, for $p\in (1, n)$, 
	\begin{eqnarray*}
		\int_{0}^{1}\frac{V([f]_t)^{\frac{np-p}{n}}}{\big[-V([f]_t)'\big]^{p-1}} \, dt&\geq& \left(\int_{0}^{1}V([f]_t)^{\frac{np-p}{n-np}}(-\,dV([f]_t) \right)^{1-p}\\
		&=& \left(\frac{np-n}{n-p}\cdot V([f]_t)^{\frac{n-p}{n-np}} \Big|_0^1\right)^{1-p}\\
		&\geq& \left(\frac{np-n}{n-p}\right)^{1-p}V([f]_1)^{\frac{n-p}{n}} \\&\geq&  \left(\frac{np-n}{n-p}\right)^{1-p}V(K)^{\frac{n-p}{n}} .
	\end{eqnarray*}  Together with (\ref{capacity-ball-1}),  Theorem  \ref{definition-smooth-1} and Corollary \ref{properties-even-affine-1}, for any $p\in (1, n)$ and for any $\tau\in [-1, 1]$, 
	\begin{eqnarray} \affp(K)&\geq& C_{p, +}(K) \nonumber \\ &=&\inf_{f\in \mathcal{F}(K)\cap C_c^{\infty}}
	\left(\int_{S^{n-1}}\lVert \nabla^{+}_u f\lVert_{p}^{-n}\ du\right)^{-\frac{p}{n}}\nonumber \\ &\geq&  n \omega_n^{\frac{p}{n}}\cdot A(n, p) \cdot \left(\frac{n-p}{p-1}\right)^{p-1}  V(K)^{\frac{n-p}{n}}\nonumber \\ &=&  A(n, p) \cdot C_p(B_n)\cdot \bigg(\frac{V(K)}{V(B_n)}\bigg)^{\frac{n-p}{n}}.\label{affine-isocapacitary-1-1} \end{eqnarray} 
	Let $K=B_n$ in inequality (\ref{affine-isocapacitary-1-1}). Then, for any $p\in (1, n)$ and  for any $\tau\in [-1,1]$, \begin{eqnarray*} \affp(B_n) \geq A(n, p) \cdot C_p(B_n). \end{eqnarray*} Together with  (\ref{inequality-affine-capacity-ball-1}), one gets, for any $p\in (1, n)$ and for any $\tau\in [-1,1]$, \begin{eqnarray}\label{affine capacity of ball} \affp(B_n) =A(n, p) \cdot C_p(B_n)= A(n, p) \cdot n\omega_n\cdot \left(\frac{n-p}{p-1}\right)^{p-1}.\label{affine-isocapacitary-2-2} \end{eqnarray} Hence, inequality (\ref{affine-isocapacitary-1-1}) can be rewritten as, for any $p\in (1, n)$, for any $\tau\in [-1, 1]$ and for any compact set $K\subset \Rn$,  $$ \bigg(\frac{\affp(K)}{\affp(B_n)}\bigg)^{\frac{1}{n-p}}\geq \bigg(\frac{V(K)}{V(B_n)}\bigg)^{\frac{1}{n}}.$$ 
	
	Now let us consider the case $p=1$. For $f\in C_c^{\infty} \cap \mathcal{F}(K)$, it can be checked, due to the dominated convergence theorem, that  for any  $u\in S^{n-1}$ and for any $\tau\in [-1, 1]$,     
	\begin{equation*}   
	\lim_{p\rightarrow 1^+} \| \varphi_\tau(\nabla_uf)\|_p=\| \varphi_\tau(\nabla_uf)\|_1.
	\end{equation*}  By Fatou's lemma, one has 
	\begin{eqnarray}    \left(\int_{S^{n-1}} \|\varphi_\tau(\nabla_uf)\|_1^{-n}\,du\right)^{-\frac{1}{n}} 
	&=&  \left(\int_{S^{n-1}}\lim_{p \rightarrow 1^+} \|\varphi_\tau(\nabla_uf)\|_p^{-n}\,du\right)^{-\frac{1}{n}} \nonumber\\ &\geq& \left(\liminf_{p \rightarrow 1^+}\int_{S^{n-1}} \|\varphi_\tau(\nabla_u f)\|_p^{-n}\,du\right)^{-\frac{1}{n}} \nonumber \\&=&  \limsup_{p \rightarrow 1^+}\left(\int_{S^{n-1}} \|\varphi_\tau(\nabla_uf)\|_p^{-n}\,du\right)^{-\frac{p}{n}}\nonumber \\ &\geq & \limsup_{p \rightarrow 1^+}\affp(K).  \nonumber
	\end{eqnarray} It follows from Theorem  \ref{definition-smooth-1}, after taking the infimum over $f \in  C_c^{\infty} \cap \mathcal{F}(K)$, that  for any $\tau\in [-1, 1]$ and for any compact set $K\subset \Rn$,    $$C_{1, \tau}(K)\geq  \limsup_{p \rightarrow 1^+} \affp(K).$$ 
	In particular, by (\ref{inequality-affine-capacity-ball-2-2}) and (\ref{affine-isocapacitary-2-2}), one has $$A(n, 1)\cdot n\omega_n\geq C_{1, \tau}(B_n)\geq  \limsup_{p \rightarrow 1^+} \affp(B_n)=A(n, 1)\cdot n\omega_n.$$ This gives the precise value of $C_{1, \tau}(B_n)$:  \begin{equation}\label{ball-affine-p-1} C_{1, \tau}(B_n)=A(n, 1)\cdot n\omega_n=\lim_{p \rightarrow 1^+} \affp(B_n),\end{equation} and hence inequality (\ref{affine-isocapacitary-1-1}) yields  $$ \bigg(\frac{C_{1, \tau}(K)}{C_{1, \tau}(B_n)}\bigg)^{\frac{1}{n-1}}\geq  \limsup_{p \rightarrow 1^+}  \bigg(\frac{\affp(K)}{\affp(B_n)}\bigg)^{\frac{1}{n-p}}  \geq \bigg(\frac{V(K)}{V(B_n)}\bigg)^{\frac{1}{n}}$$ for any $\tau\in [-1, 1]$ and for any compact set $K\subset \Rn$, as desired.  
	
	Due to the affine invariance and the translation invariance, it is trivial to see that equality holds if $K$ is an ellipsoid.\end{proof} 

 Theorem \ref{newisocapacity} asserts that the general $p$-affine capacity attains  the  minimum, among all compact sets with fixed  volume, at ellipsoids. It also asserts that ellipsoids have the maximal volumes among all compact sets with fixed general $p$-affine capacity.  When $\tau=0$, one recovers the affine isocapacitary inequality for the $p$-affine capacity proved in \cite[Theorem 3.2]{Xiao} and \cite[Theorem 1.3']{Xiao2015}. Recall that the isocapacitary inequality for the $p$-variational capacity reads: for any $p\in [1, n)$ and any compact set $K\subset \Rn$,  $$\bigg(\frac{C_p(K)}{C_p(B_n)}\bigg)^{\frac{1}{n-p}}\geq \bigg(\frac{V(K)}{V(B_n)}\bigg)^{\frac{1}{n}}. $$ It follows from Theorem \ref{comparison} that the affine isocapacitary inequality in Theorem  \ref{newisocapacity} is stronger than the isocapacitary inequality for the $p$-variational capacity. That is, for any $p\in [1, n)$, for any $\tau\in [-1, 1]$ and for any compact set $K\subset \Rn$,  $$\bigg(\frac{C_p(K)}{C_p(B_n)}\bigg)^{\frac{1}{n-p}}\geq \bigg(\frac{\affp(K)}{\affp(B_n)}\bigg)^{\frac{1}{n-p}}\geq \bigg(\frac{V(K)}{V(B_n)}\bigg)^{\frac{1}{n}}. $$ Moreover,  
combining the above inequality with \cite[(12)]{LXZ}, when  $K\subset \Rn$ is a Lipschitz star body with the origin in its interior,  the following inequality holds: for any $p\in [1, n)$ and for any $\tau\in [-1, 1]$, $$\bigg(\frac{S_p(K)}{S_p(B_n)}\bigg)^{\frac{1}{n-p}}\geq \bigg(\frac{C_p(K)}{C_p(B_n)}\bigg)^{\frac{1}{n-p}}\geq \bigg(\frac{\affp(K)}{\affp(B_n)}\bigg)^{\frac{1}{n-p}}\geq \bigg(\frac{V(K)}{V(B_n)}\bigg)^{\frac{1}{n}}, $$  where $S_p(K)$ denotes the $p$-surface area of $K$ given by formula (\ref{definition-p-surface area}).

\subsection{Connection with the general $p$-integral affine surface area} \label{section:4-1} 
In this subsection, we explore the relation between the general $p$-affine capacity and the general $p$-integral affine surface area. Throughout, denote by $\mathcal{L}_0$ the set of all Lipschitz star bodies (with respect to the origin $o$) containing $o$ in their interiors. For  a Lipschitz star body $K\in \mathcal{L}_0$, let $\nu_K(x)$ denote the unit outer normal vector of $\partial K$ at $x$ (sometimes may be abbreviated as $\nu(x)$).  Let $D_K$, the core of $K$, be given by $$D_K=\big\{ tx: \ t>0, \ x\in \partial K, \ |x\cdot \nu(x)|>0\big\}.$$ According to \cite[Lemma 5]{LXZ}, for each Lipschitz star body $K\subset \Rn$, one has $$\nu_K(x)=-\frac{\nabla \rho_K(x)}{|\nabla \rho_K(x)|} \ \ \ \mathrm{and}\ \ \ \ \nabla \rho_K(x)=-\frac{\nu_K(x)}{x\cdot \nu_K(x)}$$ for almost all $x\in \partial K\cap D_K$.

 For $p\geq 1$ and $\tau\in [-1, 1]$, define $\Pi_{p, \tau}(K)$,  the general $L_p$ projection body of $K\in \mathcal{L}_0$, to be the convex body with support function $h_{\Pi_{p, \tau}(K)}$; namely,  for any $\theta \in \sphere$,  \begin{eqnarray*} h_{\Pi_{p, \tau}(K)}(\theta)= \bigg( \int_{\partial K}\big[\varphi_\tau(\theta \cdot \nu_K(x))\big]^p\cdot |x\cdot \nu_K(x)|^{1-p}\, d\cH^{n-1}(x)\bigg)^{\frac{1}{p}}.\end{eqnarray*}
Note that $|x\cdot \nu_K(x)|^{-1}=|\nabla\rho_K(x)|$ is bounded on $\partial K$ because $\rho_K(x)$ is Lipschitz continuous on $\partial K$, and hence $h_{\Pi_{p, \tau}(K)}$ is finite.  The general $L_1$ projection body can be defined for more general sets  in $\Rn$, such as compact domains (i.e., the closure of bounded open sets) with piecewise  $C^1$ boundaries (or compact domains with finite perimeters).   When $K\in \cK_0$, formula (\ref{surfaceareameasure}) yields that, for any $\theta \in \sphere$,  \begin{eqnarray*} h_{\Pi_{p, \tau}(K)}(\theta)= \bigg( \int_{S^{n-1}}\big[ \varphi_\tau(u\cdot \theta)\big]^ph_K(u)^{1-p}\, dS(K, u)\bigg)^{\frac{1}{p}}. \end{eqnarray*}  
 Denote by $v_{p, \tau}(K, \cdot)=h^p_{\Pi_{p, \tau}(K)}(\cdot)$ the general $p$-projection function of $K$.   The general $p$-integral affine surface area of  $K\in\mathcal{L}_0$ is defined by \begin{equation}\label{definition-Phi-K}
\Phi_{p,\tau}(K)= \left(\int_{S^{n-1}}\big[v_{p,\tau}(K,u)\big]^{-\frac{n}{p}}\, du\right)^{-\frac{p}{n}}=\omega_n^{\frac{p}{n}}  V(\Pi^\circ_{p, \tau}(K)) ^{-\frac{p}{n}},
\end{equation} 
where $du$ is the normalized spherical measure and $\Pi^\circ_{p, \tau}(K)$ is the polar body of $\Pi_{p, \tau}(K)$. When $\tau=0$, one gets the $p$-integral affine surface area of $K\in \mathcal{L}_0$ in, e.g., \cite{LXZ, Zhang1}. The case $\tau=1$ defines the asymmetric $p$-integral affine surface area, denoted by $\Phi_{p,+}(K)$, of $K\in \mathcal{L}_0$. Similarly, one can also define  $\Phi_{p,-}(K)$ if $\tau=-1$. When $
K=B_n$, by (\ref{computation}), (\ref{affine capacity of ball}) and (\ref{ball-affine-p-1}), for any $p\geq 1$ and for any $\tau\in [-1, 1]$,  \begin{equation}\label{integral and capacity of ball}\Phi_{p,\tau}(B_n)=\left(\frac{n-p}{p-1}\right)^{1-p}\affp(B_n). \end{equation} It can be checked that for any $T\in GL(n)$,  $$\Phi_{p,\tau}(TK) = |\det T|^{\frac{n-p}{n}} \Phi_{p,\tau}(K).$$

Similar to the proof of Corollary \ref{properties-even-affine-1}, the following properties for the general $p$-integral affine surface area can be proved.  One cannot expect that the general $p$-integral affine surface area has the translation invariance (unless $p=1$, see following Proposition \ref{case-p=1}) and monotonicity.

\bc\label{properties-integral-affine-1} Let $p\geq 1$ and $K\in \mathcal{L}_0$. The following statements hold:  
\begin{itemize} \item [i)] for any $\tau\in [-1, 1],$ $$\Phi_{p,\tau}(K)=\Phi_{p, -\tau}(K);$$ 
	\item [ii)] for any $\lambda\in [0, 1]$ and for any $\tau, \gamma\in [-1, 1]$, $$\Phi_{p,  \lambda \tau+(1-\lambda)\gamma}(K)\geq \lambda \cdot \Phi_{p,\tau}(K)+(1-\lambda)\cdot \Phi_{p, \gamma}(K);$$
	\item [iii)]  for any $\tau\in [-1, 1]$,  $$\Phi_{p, +}(K)=\Phi_{p, -}(K)\leq \Phi_{p,\tau}(K)\leq \Phi_{p, 0}(K);$$  \item [iv)] if $-1< \tau<\gamma \leq 0,$ then $$\Phi_{p,\tau}(K)\leq\Phi_{p, \gamma}(K)$$  and  if  $0< \tau\leq \gamma< 1,$  then $$\Phi_{p, \gamma}(K)\leq \Phi_{p,\tau}(K).$$ \end{itemize}\ec 
 
By $\mathcal{C}_1$, we mean the set of all compact domains with piecewise  $C^1$ boundaries. Again, for $M\in \mathcal{C}_1$, its outer unit normal vector is denoted by $\nu_M(x)$ for $x\in \partial M$.  In the following proposition, we show that  the general $1$-affine capacity and the general $1$-integral affine surface area are all equal to the $1$-affine capacity (or equivalently, the $1$-integral affine surface area) for any $M\in \mathcal{C}_1$.

 \bp\label{case-p=1} Let $M\in \mathcal{C}_1$ be a compact domain with piecewise  $C^1$ boundary. For any $\tau\in [-1, 1]$,   one has
$$
C_{1, 0}(M)=C_{1,\tau}(M)= \Phi_{1,\tau}(M)=\Phi_{1, 0}(M). 
$$
\ep
\begin{proof} We first prove $C_{1, 0}(M)=C_{1,\tau}(M)$ for $M\in \mathcal{C}_1$;  it follows immediately from Theorem \ref{definition-smooth-1} once $\|\varphi_\tau(\nabla_u f)\|_1=\|\varphi_0(\nabla_u f)\|_1$ is established for any $f\in C_c^{\infty}\cap \mathcal{F}(M)$. To this end,  for any $M_0\in \mathcal{C}_1$ and for any $u\in \sphere$, \begin{equation} \label{centroid-1-star}
 \int_{\partial M_0}(u\cdot \nu_{M_0}(x))\,d\mathcal{H}^{n-1}(x)=0 \ \ \ \text{and} \ \ \ \int_{\partial M_0}|u\cdot \nu_{M_0}(x)|\,d\mathcal{H}^{n-1}(x)>0.\end{equation} Note that (\ref{centroid-1-star}) together with the Minkowski existence theorem leads to the powerful convexification technique, see e.g.,  \cite[p.189-190]{Zhang}.  For almost every  $t\in (0, 1)$ with  $f\in C_c^{\infty}\cap \mathcal{F}(M)$, it follows from the Sard's theorem, (\ref{Federer's coarea formula}) and  (\ref{centroid-1-star}) that, for any $\tau\in [-1, 1]$, 
\begin{eqnarray*}
\|\varphi_\tau(\nabla_u f)\|_1&=&\int_{\Rn} \bigg(\frac{1}{2}|u\cdot\nabla f|+\frac{\tau}{2}u\cdot \nabla f \bigg)\, dx\\
&=& \int_{0}^{1}\int_{\partial[f]_t}\bigg(\frac{1}{2}|u\cdot \nu(x)|+\frac{\tau}{2}u\cdot \nu(x)\bigg)\,d\mathcal{H}^{n-1}(x)dt\\
&=& \int_{0}^{1}\int_{\partial[f]_t} \frac{|u\cdot \nu(x)|}{2}\,d\mathcal{H}^{n-1}(x)dt \\&=&\int_{\Rn} \frac{|u\cdot\nabla f|}{2} \, dx\\ &=&\|\varphi_0(\nabla_u f)\|_1. 
\end{eqnarray*}   This concludes the proof of  $C_{1, 0}(M)=C_{1,\tau}(M)$ for $M\in \mathcal{C}_1$. 

  On the other hand, for any $\tau\in [-1, 1]$, 
\begin{eqnarray*} v_{1, \tau}(M, \theta)  &=&   \int_{\partial M} \varphi_\tau(\theta \cdot \nu_M(x))  \, d\cH^{n-1}(x)\\&=&   \int_{\partial M}  \bigg(\frac{|\theta \cdot \nu_M(x) |}{2}+\frac{\tau}{2} \big(\theta \cdot \nu_M(x)\big) \bigg) \, d\cH^{n-1}(x)\\ &=&  \int_{\partial M}  \frac{|\theta \cdot \nu_M(x) |}{2}  \, d\cH^{n-1}(x)\\ &=&  v_{0, \tau}(M, \theta),\end{eqnarray*} where the third equality follows again from (\ref{centroid-1-star}). Consequently, for any $\tau\in [-1, 1]$ and for any $M\in\mathcal{C}_1$,  $$\Phi_{1,\tau}(M)= \left(\int_{S^{n-1}}\big[v_{1,\tau}(K,u)\big]^{-n}\, du\right)^{-\frac{1}{n}}=\Phi_{1, 0}(M).$$

Finally,  let us  prove that $C_{1,0}(M)= \Phi_{1,0}(M)$ holds for any $M\in \mathcal{C}_1$. For each function $f\in C_c^\infty\cap \mathcal{F}(M)$, it follows from  (\ref{Federer's coarea formula}), (\ref{taufunction}), and $M\subset [f]_t$ for any $t\in [0, 1]$  that
	\begin{eqnarray*}
	\|\varphi_0(\nabla_u f) \|_1 &=&	\int_{\R^n}\varphi_0(\nabla_u f) \, dx \\ 
		&=& \frac{1}{2}\int_{0}^{1}\int_{\partial [f]_t} |u\cdot \nu(x)| \ d\mathcal{H}^{n-1}(x) \ dt\\
		&=& \frac{1}{2}\int_{0}^{1} \int_{\Pi_u [f]_t}\#([f]_t\cap(y+u\R))\,d\mathcal{H}^{n-1}(y)dt\\
		&\geq& \frac{1}{2}\int_{0}^{1} \int_{\Pi_u M}\#(M\cap(y+u\R))\,d\mathcal{H}^{n-1}(y)dt\\
		&=& \frac{1}{2}\int_{0}^{1}\int_{\partial M} |u\cdot \nu_M(x)| \ d\mathcal{H}^{n-1}(x) \ dt\\
		&=& v_{1,0}(M, u),
	\end{eqnarray*}
   where $\Pi_u K$ is the projection of $K\subset \Rn$ onto $u^{\perp}=\{x\in \Rn: \ x\cdot u=0\}$ and $\#$ denotes the number of elements of a set (see e.g., \cite{Zhang1}). Thus, for any $M\in \mathcal{C}_1$ and for any $f\in C_c^\infty\cap \mathcal{F}(M)$, 
	$$\bigg(\int_{S^{n-1}} \bigg(\int_{\R^n}\varphi_0(\nabla_uf)\ dx\bigg)^{-n}\ du\bigg)^{-\frac{1}{n}}
		\geq\bigg (\int_{S^{n-1}} v_{1,0}(M, u)^{-n}\ du \bigg)^{-\frac{1}{n}} 
		=\Phi_{1,0}(M).$$ Due to  Theorem \ref{definition-smooth-1}, by taking the infimum over $f\in C_c^{\infty}\cap \mathcal{F}(M)$, one gets, for any $M\in \mathcal{C}_1$,  
	$$
	C_{1,0}(M)\geq  \Phi_{1,0}(M).
	$$  
	
	For the opposite direction, let $\varepsilon>0$ be small enough and  consider 
	$$
	f_\varepsilon(x)=
	\begin{cases}
	0 \ \ \ \ \ \ \ \ \ \ \ \  \qquad \mathrm{if} \ \mathrm{dist}(x,K)\geq\varepsilon,\\
	1-\frac{\mathrm{dist}(x,K)}{\varepsilon} \ \ \ \ \mathrm{if} \  \mathrm{dist}(x,K)<\varepsilon. 
	\end{cases}
	$$ It has been proved in  \cite{Zhang} that  for any $u\in \sphere$, $$
	\lim_{\varepsilon \rightarrow 0}\|\varphi_0(\nabla_u f_\varepsilon)\|_1=  \lim_{\varepsilon \rightarrow 0}\int_{\R^n}\varphi_0(\nabla_u f_\varepsilon)
	\, dx=  v_{1,0}(M, u).$$ Note  that $f_\varepsilon \in \mathcal{F}(M)$ for any $\varepsilon>0$ small enough. It follows from Theorem \ref{definition2} that, for any $M\in \mathcal{C}_1$, 
	\begin{eqnarray*} 
	C_{1,0}(M)
	&\leq&  \limsup_{\varepsilon \rightarrow 0} \left(\int_{S^{n-1}} \left(\int_{\R^n}\varphi_0(\nabla_uf_\varepsilon)\ dx\right)^{-n}\ du\right)^{-\frac{1}{n}}\\ &=& \left(\liminf_{\varepsilon \rightarrow 0}  \int_{S^{n-1}} \left(\int_{\R^n}\varphi_0(\nabla_uf_\varepsilon)\ dx\right)^{-n}\ du\right)^{-\frac{1}{n}}\\ &\leq& \left(\int_{S^{n-1}} \lim_{\varepsilon \rightarrow 0}  \left(\int_{\R^n}\varphi_0(\nabla_uf_\varepsilon)\ dx\right)^{-n}\ du\right)^{-\frac{1}{n}} \\ &=& \left(\int_{S^{n-1}} \left(\lim_{\varepsilon \rightarrow 0}  \int_{\R^n}\varphi_0(\nabla_uf_\varepsilon)\ dx\right)^{-n}\ du\right)^{-\frac{1}{n}}\\ &=& \left(\int_{S^{n-1}} \left( v_{1,0}(M, u)\right)^{-n}\ du\right)^{-\frac{1}{n}}\\ &=&\Phi_{1,0}(M), 
	\end{eqnarray*} where the second inequality is due to Fatou's lemma.  This concludes the proof of $$
	C_{1,0}(M) = \Phi_{1,0}(M)
	$$  for any $M\in \mathcal{C}_1$.  \end{proof}
 When $M$ is an origin-symmetric convex body, the equality $C_{1, 0}(M)=\Phi_{1, 0}(M)$ was proved in \cite[Theorem 2]{Xiao-1}; Proposition \ref{case-p=1} extends it to all Lipschitz star bodies $M\in \mathcal{L}_0$.   The proof of $C_{1, 0}(M)=C_{1,\tau}(M)$ basically relies on the smoothness (and the convexification) of $\partial [f]_t$ instead of the compact domain $M$ itself; hence, the argument $C_{1, 0}(M)=C_{1,\tau}(M)$ holds for any compact set $M\subset \Rn$ and for any $\tau\in [-1, 1]$.  The assumption $M\in \mathcal{C}_1$ is imposed here mainly in order to have $\Phi_{1, 0}(M)$ well defined and finite. As commented in \cite[p.247]{Zhang1}, the assumption $M\in \mathcal{C}_1$  could be relaxed to more general compact domains (such as compact domains with finite perimeters). Recall the affine isoperimetric inequality for the $1$-integral affine surface area: for $M\in \mathcal{C}_1$, $$\bigg(\frac{V(M)}{V(B_n)}\bigg)^{\frac{1}{n}} \leq \bigg(\frac{\Phi_{1, 0}(M)}{\Phi_{1, 0}(B_n)}\bigg)^{\frac{1}{n-1}}$$ with equality if and only if $M$ is an ellipsoid. Then Proposition \ref{case-p=1} yields that for any $M\in \mathcal{C}_1$ and for any $\tau\in [-1, 1]$, $$\bigg(\frac{V(M)}{V(B_n)}\bigg)^{\frac{1}{n}} \leq \bigg(\frac{\Phi_{1, \tau}(M)}{\Phi_{1, \tau}(B_n)}\bigg)^{\frac{1}{n-1}}= \bigg(\frac{C_{1, \tau}(M)}{C_{1, \tau}(B_n)}\bigg)^{\frac{1}{n-1}}$$ with equality if and only if $M$ is an ellipsoid.   
 
The following theorem compares the general $p$-affine capacity and the general $p$-integral affine surface area. We only concentrate on $p\in (1, n)$ as the case $p=1$ has been discussed in Proposition \ref{case-p=1}. When $\tau=0$ and $K$ is an origin-symmetric convex body, it recovers \cite[Theorem 3.5]{Xiao}. 
\bt \label{capacityandintegral}
Let  $K\in \mathcal{L}_0$ and  $1< p<n$. The following inequality
$$
\frac{C_{p,\tau}(K)}{C_{p,\tau}(B_n)}\leq \frac{\Phi_{p,\tau}(K)}{\Phi_{p,\tau}(B_n)} \ \ \ \ \mathrm{for \ \ any} \ \ \tau\in [-1, 1]
$$ holds 
with equality if $K$ is an origin-symmetric ellipsoid.
\et  
\begin{proof} Let $K\in \mathcal{L}_0$ and $p\in (1, n)$. Define the function $g$ by: for $s>0$, $$g(s)=\mathrm{min}\Big\{1,\ \  s^{\frac{n-p}{1-p}}\Big\}.$$ Let $f(x)=g\big(\frac{1}{\rho_K(x)}\big). $ Then  $f(x)\geq \textbf{1}_K$  and  $\|f\|_{\infty}=1$. From (\ref{levelset}) and the fact that $g$ is strictly decreasing on $s\in (1, \infty)$, it follows that, for all $t\in (0, 1)$ with $t=g(s)=s^{\frac{n-p}{1-p}}$,  $$[f]_t=\{x\in \R^n: 1/\rho_K(x)\leq s\}. $$ That is, $[f]_t=[f]_{g(s)}=sK$ for any $s>1$. Together with \cite[Lemma 6]{LXZ},  for any $x\in \partial[f]_t$, there exists $z\in \partial K$ with $x=sz$ such that   $$|\nabla f(x)|=\frac{|g'(s)|}{|z\cdot \nu_K(z)|} \ \ \ \ \mathrm{and} \ \ \  \ \nu_K(z)=\nu_{[f]_t} (x)=-\frac{\nabla f(x)}{|\nabla f(x)|}.
$$
 By (\ref{Federer's coarea formula}), one has,  for any $u\in \sphere$, \begin{eqnarray*}
		\| \varphi_\tau(\nabla_u f)\|_p^p
		&=& \int_{0}^{1}\int_{\partial [f]_t} \big[ \varphi_\tau(-u\cdot \nu_{[f]_t} (x))\big]^p\cdot |\nabla f(x)|^{p-1}\, d\mathcal{H}^{n-1}(x)\, dt \\ &=& \int_{1}^{\infty}|g'(s)|\int_{\partial [f]_{g(s)}} \big[\varphi_\tau(-u\cdot\nu_{[f]_{g(s)}} (x))\big]^p\cdot |\nabla f(x)|^{p-1}\, d\mathcal{H}^{n-1}(x)\, ds\\ &=& \int_{1}^{\infty} |g'(s)|^{p} s^{n-1}  \int_{\partial K}  \big[\varphi_\tau(-u\cdot \nu_K(z))\big]^p \cdot |z\cdot \nu_K(z)|^{1-p}\, d\mathcal{H}^{n-1}(z)\,ds \\ &=& \bigg(\frac{n-p}{p-1}\bigg)^{p} \bigg(\int_{1}^{\infty}  s^{\frac{n-1}{1-p}} \,ds \bigg) \bigg( \int_{\partial K}  \big[ \varphi_\tau(-u\cdot \nu_K(z))\big]^p \cdot |z\cdot \nu_K(z)|^{1-p}\, d\mathcal{H}^{n-1}(z)\bigg) \\&=&  \bigg(\frac{n-p}{p-1}\bigg)^{p-1}v_{p, \tau}(K, -u).\end{eqnarray*} It follows from  (\ref{definition-h-f}) and (\ref{definition-Phi-K})  that \begin{eqnarray*}\mathcal{H}_{p,\tau}(f)&=&
		 \bigg(\int_{S^{n-1}}\lVert \varphi_\tau(\nabla_u f)\lVert_p^{-n}\ du\bigg)^{-\frac{p}{n}} \\
		&=&\left(\frac{n-p}{p-1}\right)^{p-1}\bigg(\int_{S^{n-1}} v_{p, \tau}(K, -u)^{-\frac{n}{p}}\ du\bigg)^{-\frac{p}{n}}\\
		&=& \left(\frac{n-p}{p-1}\right)^{p-1}\Phi_{p,\tau}(K).
	\end{eqnarray*}	 A standard limiting argument together with Definition \ref{definition1}  show that, for any $p\in (1, n)$, for any $\tau\in [-1, 1]$ and for any $K\in \mathcal{L}_0$, 
	\begin{eqnarray*}
		C_{p,\tau}(K) \leq   \left(\frac{n-p}{p-1}\right)^{p-1}\Phi_{p,\tau}(K).
	\end{eqnarray*}
	By (\ref{integral and capacity of ball}), the above inequality can be rewritten as
	$$
	\frac{C_{p,\tau}(K)}{C_{p,\tau}(B_n)}\leq \frac{\Phi_{p,\tau}(K)}{\Phi_{p,\tau}(B_n)}.
	$$

	Clearly equality holds in the above inequality if $K=B_n$. Due to the affine invariance of both $\affp(\cdot)$ and $\Phi_{p, \tau}(\cdot)$, equality holds in the above inequality if $K$ is an origin-symmetric ellipsoid.  
\end{proof}

  Together with \cite[(13)]{LXZ}, Corollary \ref{properties-integral-affine-1}  and Theorem \ref{newisocapacity}, for any $K\in \mathcal{L}_0$, for any $p\in (1, n)$ and for any $\tau\in [-1, 1]$,  one has, \begin{equation}\label{zhang-affine-iso-1} \bigg(\frac{S_{p}(K)}{S_{p}(B_n)}\bigg)^{\frac{1}{n-p}}\geq\bigg(\frac{\Phi_{p, 0}(K)}{\Phi_{p, 0}(B_n)}\bigg)^{\frac{1}{n-p}}\geq \bigg(\frac{\Phi_{p,\tau}(K)}{\Phi_{p,\tau}(B_n)}\bigg)^{\frac{1}{n-p}}\geq \bigg(\frac{\affp(K)}{\affp(B_n)}\bigg)^{\frac{1}{n-p}}\geq \bigg(\frac{V(K)}{V(B_n)}\bigg)^{\frac{1}{n}} \end{equation}  
with equality if $K$ is an origin-symmetric ellipsoid. Inequality (\ref{zhang-affine-iso-1})  extends several known results in the literature. For example, inequality (\ref{zhang-affine-iso-1})  strengthens the following (affine) isoperimetric inequality  (see \cite[inequality (13)]{LXZ}): for $\tau=0$ and for any $K\in \mathcal{L}_0$, 
$$\bigg(\frac{S_p(K)}{S_p(B_n)}\bigg)^{\frac{1}{n-p}}\geq   \bigg(\frac{\Phi_{p,0}(K)}{\Phi_{p, 0}(B_n)}\bigg)^{\frac{1}{n-p}}\geq  \bigg(\frac{V(K)}{V(B_n)}\bigg)^{\frac{1}{n}}.  $$   Moreover, inequality (\ref{zhang-affine-iso-1})  holds for all $K\in \cK_0\subset \mathcal{L}_0$, and hence it   extends the following affine isoperimetric inequality (\ref{integral-isoperimetric-inequality}) for convex bodies  to Lipschitz star bodies: for any  $K\in \cK_0$, for any $\tau\in [-1, 1]$ and for any $p\in (1, n)$,   \begin{equation}\label{integral-isoperimetric-inequality} \bigg(\frac{\Phi_{p,\tau}(K)}{\Phi_{p,\tau}(B_n)}\bigg)^{\frac{1}{n-p}}\geq \bigg(\frac{V(K)}{V(B_n)}\bigg)^{\frac{1}{n}},\end{equation} which is an immediate consequence of the general $L_p$ affine isoperimetric inequality for the general $L_p$ projection body \cite{HSgeneral}.

\vskip 2mm   \noindent {\bf Acknowledgments.}  DY is supported by a NSERC
grant.

\vskip 2mm \noindent Han Hong, \ \ \ {\small \tt honghan0927@gmail.com}\\
{ \em Department of Mathematics and Statistics,   Memorial University of Newfoundland,
   St.\ John's, Newfoundland, Canada A1C 5S7 }

\vskip 2mm \noindent Deping Ye, \ \ \ {\small \tt deping.ye@mun.ca}\\
{ \em Department of Mathematics and Statistics,
   Memorial University of Newfoundland,
   St.\ John's, Newfoundland, Canada A1C 5S7 }


\begin{thebibliography}{99}
 
\bibitem{Borell}
C. Borell, {\em Capacitary inequalities of the Brunn-Minkowski type}, Math. Ann. 263 (1983) 179-184.

\bibitem{Caffarelli1996} L.A. Caffarelli, D. Jerison and E.H. Lieb,  {\em On the case of equality in the Brunn-Minkowski inequality for capacity,}  Adv. Math. 117 (1996) 193-207.

\bibitem{CLYZ} A. Cianchi, E. Lutwak, D.Yang and G. Zhang, {\em Affine Moser-Trudinger and Morrey-Sobolev inequalities,} Calc. Var. Partial Differ. Equ. 36 (2009) 419-436.

\bibitem {CNSXYZ} A. Colesanti, K. Nystr\"om, P. Salani, J. Xiao, D. Yang and G. Zhang, {\em The Hadamard variational formula and the Minkowski problem for $p$-capacity,} Adv. Math. 285 (2015) 1511-1588.


\bibitem{CS2003} A. Colesanti and P. Salani, {\em The Brunn–-Minkowski inequality for $p$-capacity of convex bodies}, Math. Ann. 327 (2003) 459–-479.

\bibitem{Evan} L. Evans and R. Gariepy, {\em Measure Theory and Fine Properties of Functions,} CRC Press LLC. 1992.

\bibitem{Federer} H. Federer, {\em Geometric Measure Theory,} Springer, Berlin, 1969.


\bibitem{Gardner-1} R.J. Gardner, {\em The Brunn-Minkowski inequality,} Bull. Amer. Math. Soc. 39 (2002) 355-405.

\bibitem{GHW2014} R.J. Gardner, D. Hug and W. Weil, {\em The Orlicz-Brunn-Minkowski theory: a general framework, additions, and inequalities,}  J. Differential Geom.  97 (2014) 427-476.

\bibitem{ghwy15} R.J. Gardner, D. Hug, W. Weil and D. Ye, {\em The dual
Orlicz-Brunn-Minkowski theory,} J. Math. Anal. Appl. 430 (2015) 810-829.


\bibitem{Groemer} H. Groemer, {\em Geometric applications of Fourier series and spherical harmonics,} Cambridge University Press, New York, 1996.

\bibitem{HSgeneral} C. Haberl and F.  Schuster, {\em General $L_p$ affine isoperimetric inequalities,} J. Differential. Geom. 83 (2009) 1-26.

\bibitem{HSasymmetric} C. Haberl and F. Schuster, {\em Asymmetic affine $L_p$ Sobolev inequalities,} J. Funct. Anal. 257 (2009) 641-658.

\bibitem{HSXasymmetric} C. Haberl, F. Schuster and J. Xiao, {\em An asymmetic affine P\'{o}lya-Szeg\"{o} Principle,} Math. Ann. 352 (2012) 517-542.

\bibitem{HJM2016}  J. Haddad, C.H. Jim\'{e}nez and M. Montenegro, {\em 
 Sharp affine Sobolev type inequalities via the $L_p$ Busemann-Petty centroid inequality,} J. Funct. Anal. 271 (2016) 454-473. 

\bibitem{HKM-1} J. Heinonen, T. Kilpel\"{a}inen and O. Martio, {\em Nonlinear Potential Theory of Degenerate Elliptic Equations,} Dover Publications, 2006. 

 
\bibitem{HYZ-17} H. Hong, D. Ye and N. Zhang, {\em The $p$-capacitary Orlicz-Hadamard variational formula and Orlicz-Minkowski problems,} submitted. 
 

\bibitem{Jerison}
  D. Jerison, {\em A Minkowski problem for electrostatic capacity,} Acta Math. 176 (1996) 1-47.

\bibitem{Jerison-1996}
  D. Jerison, {\em The Direct Method in the Calculus of Variations for Convex Bodies}, Adv. Math. 122 (1996)  262-279. 
 
 
\bibitem{Ldvaluation} M. Ludwig, {\em Minkowski valuations,} Trans. Amer. Math. Soc. 357 (2005) 4191-4213.

\bibitem{Ludwig2010} {M. Ludwig,} {\em General affine surface areas}, Adv. Math. 224
(2010)  2346-2360.

\bibitem{LXZ} M. Ludwig, J. Xiao and G. Zhang, {\em Sharp convex Lorentz-Sobolev inequalities,} Math. Ann. 350 (2011) 169-197.

\bibitem{Lutwak} E. Lutwak, {\em The Brunn-Minkowski-Fiery theory. I. Mixed volumes and the Minkowski problem,} J. Differential. Geom. 38 (1993) 131-150.

\bibitem{LYZ1} E. Lutwak, D. Yang and G. Zhang, {\em $L_p$ affine isoperimetric inequalities,} J. Differential. Geom. 56 (2000) 111-132.

\bibitem{LYZ-duke} E. Lutwak, D. Yang and G. Zhang, {\em A new ellipsoid associated with convex bodies,} Duke Math. J. 104 (2000) 375-390. 


\bibitem{LYZ} E. Lutwak, D. Yang and G. Zhang, {\em Sharp affine $L_p$ Sobolev inequalites,} J. Differential. Geom. 62 (2002) 17-38.

\bibitem{LYZ2010a}
 E. Lutwak, D. Yang and G. Zhang, {\em Orlicz projection bodies,}
  Adv. Math. {223} (2010)  220-242.
  
  \bibitem {LYZ2010b}  E. Lutwak, D. Yang and G. Zhang, {\em Orlicz centroid bodies,}  J. Differential. Geom. 84  (2010)  365-387. 
  
 
\bibitem{Mazya-85} V. Maz'ya, {\em Sobolev spaces,} Springer-Verlag, Berlin,  1985. 


\bibitem{Mazya} V. Maz'ya, {\em Sobolev Spaces with Applications to Elliptic Partial Differential Equations,} 2nd, Springer, 2011.


\bibitem{VHN} V.H. Nguyen, {\em New approach to the affine P\'{o}lya-Szeg\"{o} principle and the stability version of the affine Sobolev inequality,} Adv. Math. 302 (2016) 1080-1110.

\bibitem{Ober} M. Ober, {\em Asymmetric $L_p$ convexification and the convex Lorentz Sobolev inequality,} Monatsh. Math. 179 (2016) 113-127.

\bibitem{Petty1971}
C. Petty,  {\em Isoperimetric problems,} Proc. Conf. Convexity and Combinatorial
Geometry, Univ. of Oklahoma, Norman, 1971, 26-41.

  
\bibitem{Schneider} R.\ Schneider, {\em Convex bodies: The Brunn-Minkowski theory,} Cambridge Univ. Press, Cambridge, 1993.
 
\bibitem{WangTuo2012} T. Wang, {\em The affine Sobolev-Zhang inequality on $BV(\Rn)$,} Adv. Math. 230 (2012) 2457-2473.



\bibitem{WangTuo2013} T. Wang,  {\em The affine P\'{o}lya-Szeg\"{o} principle:
Equality cases and stability,}    J. Funct. Anal. 265 (2013) 1728-1748.

\bibitem{WangTuo2015} T. Wang,  {\em On the Discrete Functional $L_p$ Minkowski Problem,}  Int. Math. Res. Notices, 2015  (2015)  10563-10585. 


\bibitem{WangTuo2017} T. Wang and J. Xiao, {\em The affine BV-capacity,} {\tt arXiv:1510.07920.}  

\bibitem{Weberndorfer} M. Weberndorfer, {\em Shadow systems of asymmetric $L_p$ zonotopes,} Adv. Math. 240 (2013) 613-635.

\bibitem{XJL} D. Xi, H. Jin  and G. Leng, {\em The Orlicz Brunn-Minkowski inequality},  Adv. Math.  260 (2014) 350-374.
 
 
\bibitem{Xiao2007} J. Xiao, {\em The affine Sobolev and isoperimetric inequalities split twice,} Adv. Math. 211 (2007) 417-435. 

\bibitem{Xiao2015} J. Xiao, {\em Corrigendum to ``The sharp Sobolev and isoperimetric inequalities split twice"}, Adv. Math. 268 (2015) 906-914.

\bibitem{Xiao}	J. Xiao, {\em The p-Affine capacity,} J. Geom. Anal. 26 (2016) 947-966.

\bibitem{Xiao-1} J. Xiao, {\em The p-Affine Capacity Redux,}  J. Geom. Anal.,  in press. 

\bibitem{Xiaozhang} J. Xiao and N. Zhang, {\em The Relative p-Affine capacity,} Proc. Amer. Math. Soc. 144 (2016) 3537-–3554.

  
\bibitem{Zhang} G. Zhang, {\em The affine Sobolev inequality,} J. Differential. Geom. 53 (1999) 183-202.

\bibitem{Zhang1} G. Zhang, {\em New affine isoperimetric inequalities,} ICCM. Vol. II (2007) 239-267.

\bibitem{Zhub2014} B. Zhu, J. Zhou and W. Xu, {\em Dual Orlicz-Brunn-Minkowski theory},
Adv.  Math. 264 (2014) 700-725.

\bibitem{xiongzou} D. Zou and G. Xiong,  {\em A unified treatment for $L_p$ Brunn-Minkowski type inequalities,} Commun. Anal. Geom., in press. 

\end{thebibliography}
\end{document}